\documentclass[11pt,a4paper]{article}
\usepackage[top=0.75 in,bottom=0.75 in, left=0.75 in, right=0.75 in]{geometry}   
\usepackage{bm}
\usepackage{graphicx}
\usepackage[font=small]{caption}
\usepackage{amssymb}
\usepackage{amsthm}
\usepackage{amsmath}
\usepackage{amssymb}
\usepackage{stmaryrd}
\usepackage{epstopdf}
\usepackage{subcaption}
\usepackage{sidecap}
\usepackage{authblk}
\usepackage{enumitem}
\usepackage{color}
\usepackage{ulem}
\usepackage{comment}
\usepackage{float}
\usepackage{soul}
\usepackage{hyperref}
\usepackage{enumitem}

\theoremstyle{plain}
\newtheorem{theorem}{Theorem}[section]
\newtheorem{corollary}{Corollary}[section]

\newtheorem{proposition}{Proposition}[section]
\theoremstyle{definition}

\theoremstyle{remark}

\DeclareMathOperator{\vol}{len}

\DeclareMathOperator{\supp}{supp}

\title{Irreversibility condition and stability of equilibria in the inverse-deformation approach to fracture} 
\author[1,*]{Arnav Gupta}
\affil[1]{\small{Field of Theoretical and Applied Mechanics, Cornell University, Ithaca, NY 14853, USA}}
\affil[*]{\small{Corresponding author e-mail: ag2592@cornell.edu}}
\date{}

\begin{document}
\maketitle
\begin{center}\vspace{-2em}
	\textit{Dedicated to Timothy J. Healey on the occasion of his 70\textsuperscript{th} birthday}
\end{center}
\begin{abstract}
We derive an irreversibility condition for fracture in the inverse-deformation approach using the second law of thermodynamics. We consider the problem of brittle failure in one-dimensional, homogeneous, elastic bars undergoing quasistatic motion. Despite the presence of a non-zero interfacial/surface energy, the third derivative of the inverse deformation is discontinuous at crack faces. This is due to the inequality constraint allowing fracture while ensuring that the orientation of matter is preserved. We show that any change in the location of open cracks in their reference configuration result in negative entropy production. We disallow these changes according to the second law of thermodynamics, fixing the location of cracks in their reference domain. The inequality constraint and the irreversibility condition limit the space of admissible variations. We prove second-order conditions for local stability incorporating these restrictions. Their numerical implementation shows that all broken equilibria of the elastic bar under hard loading are locally stable.  
\end{abstract}

\numberwithin{equation}{section}  
\renewcommand{\vec}[1]{\bm{#1}}

\tableofcontents
\newcommand{\sw}{{W^\ast}}
\newcommand{\ssw}{{\buildrel \ast\ast\over W}}
\newcommand{\beq}{\begin{equation}}
\newcommand{\eeq}{\end{equation}}
\newcommand{\la}{{\lambda}}
\newcommand{\Om}{{\Omega}}
\newcommand{\Oma}{{\Omega_\ast}}
\newcommand{\Oms}{{\Omega_s}}
\newcommand{\Omn}{{\Omega_n}}

\newcommand{\wom}{{\widehat{\Omega}}}
\newcommand{\woma}{{\widehat{\Omega}_\ast}}

\newcommand{\ppwom}{\partial_1{\widehat{\Omega}}}

\newcommand{\pppwom}{\partial_2{\widehat{\Omega}}}

\newcommand{\pwom}{\partial{\widehat{\Omega}}}

\newcommand{\po}{{\partial\Omega}}
\newcommand{\poa}{{\partial\Omega_\ast}}
\newcommand{\pd}{\partial}
\newcommand{\si}{{\sigma}}
\newcommand{\TT}{\mathcal{T}}
\newcommand{\ep}{{\varepsilon}}

\newcommand{\epo}{{\varepsilon_0}}

\newcommand{\xwom}{{\chi_{{}_\wom}}}

\newcommand{\epi}{{\varepsilon_1}}
\newcommand{\bpr }{\begin{proof}}
\newcommand{\epr}{\end{proof}}
\newcommand{\tred}{\textcolor{black}}
\newcommand{\vp}{{\varpi}}
\newcommand{\bF}{\bm{F}}
\newcommand{\bA}{\bm{A}}
\newcommand{\bS}{\bm{S}}
\newcommand{\bP}{\bm{P}}
\newcommand{\bT}{\bm{T}}
\newcommand{\bG}{\bm{G}}
\newcommand{\bK}{\bm{K}}
\newcommand{\bH}{\bm{H}}
\newcommand{\bff}{\bm{f}}
\newcommand{\bh}{\bm{h}}
\newcommand{\bw}{\bm{w}}
\newcommand{\bu}{\bm{u}}
\newcommand{\bx}{\bm{x}}
\newcommand{\by}{\bm{y}}
\newcommand{\bi}{\bm{I}}
\newcommand{\bg}{\bm{g}}
\newcommand{\bz}{\bm{z}}
\newcommand{\bs}{\bm{s}}
\newcommand{\id}{\bm{1}}

\newcommand{\cD}{{\mathcal{D}}}

\newcommand{\cDa}{{\mathcal{D}_\ast}}

\newcommand{\re}{{\mathbb{R}}}
\newcommand{\rt}{{\mathbb{R}^2}}
\newcommand{\lp}{{\mathrm{Lin}_+}}
\newcommand{\lpo}{{\mathrm{Lin}^0_+}}

\newcommand{\lo}{{\mathrm{Lin}^0}}

\newcommand{\lin}{{\mathrm{Lin}}}

\section{Introduction}\label{intro}
We derive an irreversibility condition for fracture in the inverse-deformation approach \cite{RHA} using the second law of thermodynamics. We consider a one-dimensional, homogeneous, elastic bar undergoing quasistatic motion. The inverse deformation maps the deformed configuration of the bar to the reference domain. Despite the presence of a non-zero interfacial/surface energy, the third derivative of the inverse deformation is discontinuous at crack faces. This is due to the inequality constraint ensuring that the inverse strain\footnote{We use the term ``strain" to denote the deformation gradient in this one-dimensional setting. The term ``inverse strain" is the gradient of the inverse deformation.} is nonnegative and the orientation of matter is preserved. We show that any change in the location of open cracks in their reference configuration result in negative entropy production. We disallow these changes according to the second law of thermodynamics, fixing the location of cracks in their reference domain.

The inverse-deformation approach to fracture for one-dimensional, elastic bars is proposed in \cite{RHA}. In that work, the inverse deformation is extended to a continuous, piecewise smooth function even in the presence of cracks. A local, nonconvex constitutive model is employed in \cite{RHA} and regularized using the gradient of inverse strain. This equips the model with a surface energy that also embodies Griffith's constant energy in a certain Gamma-limit \cite{HPP}. A global bifurcation analysis in \cite{RHA} identifies all equilibria of the elastic bar under hard loading. The bar deforms homogeneously until it fractures due to a snap instability at a critical load. Spontaneous nucleation of sharp cracks is predicted in \cite{RHA} without using any damage/phase variables or pre-existing flaws. The stress in the bar drops exactly to zero and the bar is broken. The crack faces are sharp and the actual deformation is discontinuous. The present work is a follow-up to \cite{RHA} that establishes the nontriviality of the second law of thermodynamics in this method.

The analyses in \cite{RHA} predict that nucleation of a crack occurs at one of the two ends of the bar. They report that all other broken equilibria, that have internal or multiple cracks, are unstable and unobservable in reality. We address this question in this article. We incorporate the irreversibility condition derived here in the stability analyses of broken equilibria. We prove second-order conditions for local stability forbidding the variations that change the location of cracks. Their numerical implementation shows that all broken equilibria of the elastic bar under hard loading are locally stable. The elastic bar can lock into any one of these solutions, resulting in internal or multiple cracks upon failure.

\subsection*{Irreversibility condition}
Irreversibility is introduced as an additional postulate in Griffith's theory \cite{griffith1921, griffith1924theory}. In modern terms, it is defined as follows: If at time $t$, $\Gamma(t)$ is the set of all material points on a crack surface, then
	$$\Gamma(s)\subseteq \Gamma(t) \text{ for all } s\leq t.$$ 
Various equilibrium cracks under loading are unstable in the absence of such a condition; cracks cannot exist in the Griffith's theory without postulating a mechanism for irreversibility\footnote{See section E in \cite{rice1968}.}. In the Linear Elastic Fracture Mechanics (LEFM), irreversibility is incorporated by postulating pre-existing cracks that never shorten in length. 

The same postulate is extended to finite-deformation, nonlinear elasticity with pre-existing cracks in \cite{le, stumpf}. Irreversibility is assumed in these works by restricting the admissible space. A broken configurations is admissible if the new surface of discontinuity contains or coincides with that of all previous states. Similar irreversibility condition is also assumed in \cite{gurtin2000, gurtin1996, gurtin1998}.

Irreversibility is assumed in the phase-field methods in the form of an evolution law in \cite{francmarigo98, bourdin}. In \cite{bourdin08}, it is assumed that the `Griffith-like' potential is infinity if the crack set $\Gamma(t)$ decreases. In \cite{fried2013}, the change in phase variable is \textit{a priori} constrained to be unidirectional.

Irreversibility in the inverse-deformation approach to fracture arises as a consequence of the second law of thermodynamics. We require no additional assumption to prove this condition. For a broken, one-dimensional bar, the inverse deformation and its first two derivatives are continuous while its third derivative is discontinuous at the crack faces. The entropy of the system changes if such discontinuities move in the reference configuration.

Entropy production due to the motion of discontinuities at interfaces is studied in the context of phase transitions in \cite{abeyknowles88, abeyknowlesbook}. In the absence of interfacial energy, a multiphase solution has discontinuities in strain. The motion of such discontinuities in the reference configuration causes change in the entropy of the system. We do a similar analysis here for fracture in the inverse-deformation formulation.

We summarize the inverse-deformation approach to fracture from \cite{RHA} in section \ref{invdef}. Following the same work, we incorporate a surface energy in the model that is quadratic in inverse-strain gradient. In phase-transition models, the inclusion of interfacial energy in terms of strain gradients regularizes the strain discontinuities and results in classical equilibrium solutions \cite{carr}. This is not the case for a broken solution here: We prove in section \ref{smooth} that the third derivative of the inverse-deformation is discontinuous. 

In section\ \ref{workext}, we derive an expression for work done on a part of the bar by forces external to it. We employ this in section \ref{entropy} while calculating entropy production. We consider an elastic bar undergoing quasistatic motion. Inertia affects are not included and the notion of time plays the role of a history parameter \cite{knowles79}. The entropy production near the crack face is proportional to the jump in the third derivative times the velocity of the crack face in the reference domain. The jumps in the third derivative are related with the positive Lagrange-multiplier field enforcing the inequality constraint and their sign are easily determined. We show that this velocity of the crack in the reference domain necessarily equals zero to ensure that the second law of thermodynamics is satisfied. 

The irreversibility condition derived here does not preclude healing. The entropy calculations in section \ref{entropy} require that the crack faces do not touch in the deformed configuration, i.e.\ the cracks are open. The behavior of crack faces in contact depends on other factors such as their atomistic structure and surrounding environment. For example, in the presence of air or humidity the newly formed surfaces of a cracked metal bar oxidize and are difficult to adhere. In contrast, atomically clean surfaces in ultrahigh-vacuum join with each other due to a process called cold welding \cite{lu2010cold}. The reversibility of fracture in mica when crack faces are in contact is shown in \cite{obreimoff}.

Furthermore, irreversibility does not restrict the location of a crack face in the deformed configuration. Their spatial location may change but the inverse-deformation consistently maps the crack face to the same breakage point in the reference configuration. 

\subsection*{Stability of broken equilibria}
Stability in the absence of irreversibility was first studied in \cite[Proposition 4.7]{RHA}. The irreversibility condition was not known to the authors. In its absence, they show that solutions with a crack at one end of the bar are the only locally stable, broken solutions. Using the method in \cite[Section 8]{carr}, they construct specific variations of the inverse deformation for which the second variation of the total energy is negative for all other broken equilibria. In the phase-transition problem studied in \cite{carr}, the first variation of the energy at an equilibrium is zero for all admissible variations. This is not true in the fracture problem studied in \cite{RHA} due to the presence of the inequality constraint. For solutions with two end cracks, the first variation of the energy is positive for the variations constructed in \cite[Proposition 4.7]{RHA}. This renders second-order calculations in those directions inconsequential to the stability of these equilibria. 

We refer to the stability criterion in the absence of irreversibility as Criterion 1. We prove a second-order condition for local stability subject to this criterion in section\ \ref{c1}. This condition incorporates the possibility of positive first variation of the energy at a broken equilibrium. We use this condition to modify \cite[Proposition 4.7]{RHA} in Proposition\ \ref{c1p1} and show that any broken solution with an internal crack is unstable. The solutions with cracks at either one end of the bar or both ends of the bar remain. The stability of these solution is the same in the presence or absence of irreversibility. In section\ \ref{stabresults}, we find that these solutions with one or two end cracks are both locally stable in the absence of the irreversibility condition. 
 
Notwithstanding that, the purpose of this work is to determine the stability of solutions in the presence of the irreversibility condition. We refer to this criterion as Criterion 2 and study it in section\ \ref{c2}. We show that an equilibrium is locally stable if the second variation of the total energy is positive for all admissible variations that are zero in the empty regions. These empty regions are the spaces in the deformed configuration created due to the opening of cracks. This result means, for the purposes of stability, every connected, unbroken part of the elastic bar is decoupled. The full broken solution is locally stable if such unbroken parts of the bar are all \textit{independently} stable. An unbroken part is \textit{independently} stable if the energy increases locally for all admissible variation that are non-zero only within that part. 

An unbroken part with crack faces on both ends can also freely translate within the empty region. At a given stretch, these translations result into a one-parameter family of equilibria that have the same total energy. If the second variations of the energy is strictly positive for variations that are devoid of such translations, the unbroken part is stable. 
 
The variations constructed in section\ \ref{c1}, that render solutions with internal cracks unstable in the absence of irreversibility, change the location of these cracks. These variations are inadmissible in Criterion 2. Furthermore, the stability analysis using Criterion 2 satisfies the second law of thermodynamics and models the real-world problem. We implement the stability condition due to Criterion 2 using the finite-element discretization of \cite{GH} in section\ \ref{num}. We find that all broken solutions of the elastic bar are locally stable. The global minima appears in the first bifurcating branch with one end crack and is energetically preferred. This is also shown in the Gamma-limit of the small parameter characterizing the surface energy approaching zero in \cite{HPP}. The remaining broken solutions, with internal or multiple cracks, appear in the higher modes and are local minimizers.

\section{Preliminaries}\label{prelimin}

\subsection{Inverse-deformation formulation} \label{invdef}

\begin{figure}[h]
	\centering
	\begin{subfigure}[c]{0.4\textwidth}
		\centering
         \includegraphics[width=\textwidth]{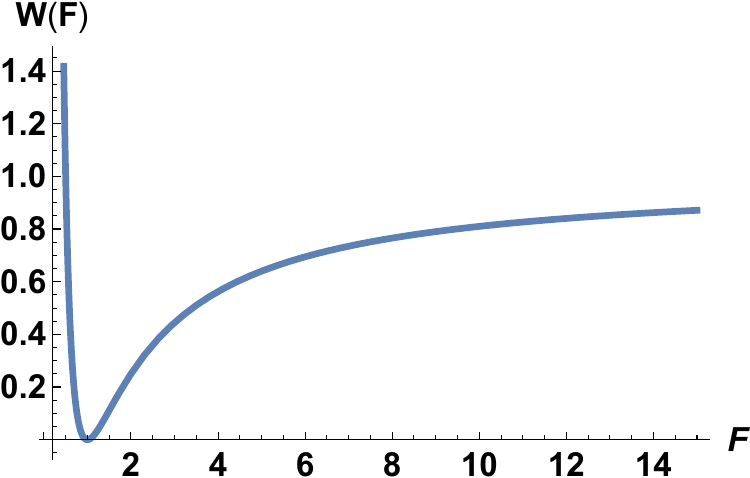}
         \caption{}
          \label{W}
	\end{subfigure}
	\begin{subfigure}[c]{0.4\textwidth}
		\centering
         \includegraphics[width=\textwidth]{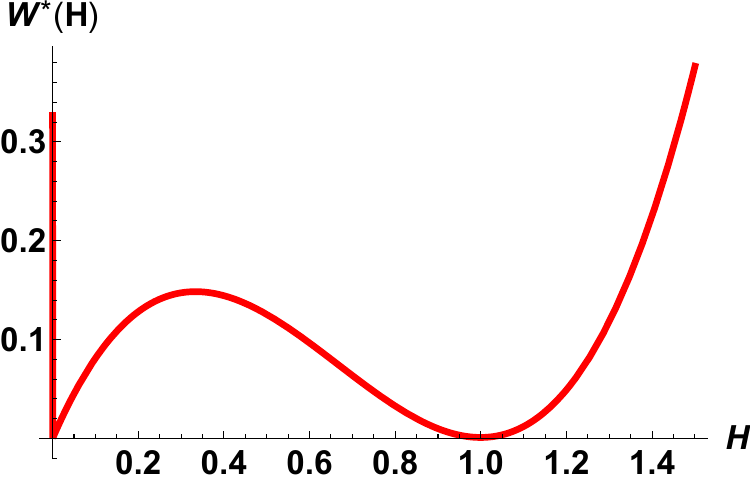}
         \caption{}
          \label{Wstar}
	\end{subfigure}
  \caption{(a) The stored-energy density function $W$ for a material that undergoes brittle fracture. (b) The corresponding inverse-energy density function $W^*$ with wells at $H=0$ and $H=1$. These figures are reproduced here from \cite{RHA} with permission.}
  \label{energydens}
\end{figure}
  
Consider a bar occupying the region $\Om:=[0,1]$ in the reference configuration. We fix one end of the bar and stretch the other to a length $\lambda$ using a hard-loading device. The deformation $f: \Om\rightarrow \mathbb{R}$ satisfies the boundary conditions 
\begin{equation}
	f(0)=0,\quad f(1)=\lambda.
\end{equation}
The stored-energy density function $W:(0,\infty)\rightarrow [0,\infty)$ of a brittle material is a Lennard-Jones type potential as shown in Figure\,\ref{W}. It asymptotes to a constant $\gamma$ as the actual strain approaches infinity. It is assumed as in \cite{RHA} that 
\begin{equation}
	\begin{aligned} \label{hypo}
	& W\in C^2(0,\infty);\, W(F)\nearrow \infty \text{ as } F\searrow 0;\, W(F)\nearrow\gamma>0 \text{ as } F\nearrow\infty;\\
	& W(1)=0;\, W(F)>0 \text{ for all } F\neq 1.
\end{aligned}
\end{equation}
The total energy of the body is
\begin{equation}\label{E}
	E[\lambda,f]=\int_0^1 W(f'(x))\,dx.
\end{equation}

\begin{figure}[h]
  \centering
\includegraphics[width=0.85\textwidth]{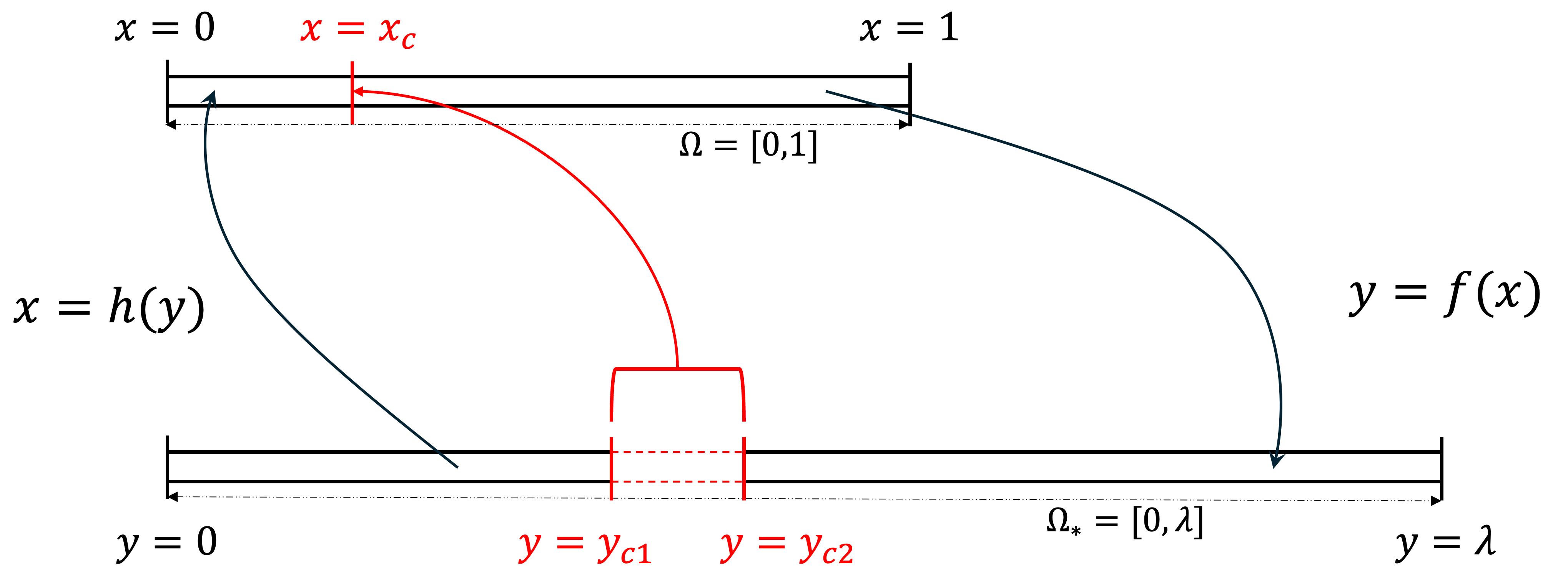}
  \caption{The deformation $f$ maps  the reference $\Om$ onto the deformed configuration $\Oma$. The inverse deformation $h$ maps the deformed configuration $\Oma$ onto the reference domain $\Om$. The region $[y_{c1},y_{c2}]\subset \Oma$ shows empty space due to opening of a crack. The inverse deformation $h$ maps all points in $[y_{c1},y_{c2}]$ to the crack location $x=x_c$ in $\Om$.}
  \label{inverse}
\end{figure}

The deformed configuration is defined as $\Oma:=[0,\lambda]$. Let $h:\Oma\rightarrow \Om$ be the inverse deformation, as shown in Figure \ref{inverse}. We use the change of variable $x=h(y),\, f'=1/h'$ in \eqref{E} following \cite{shield} and redefine the total energy
\begin{equation}\label{Es}
	E^*[\lambda,h]:=E[\lambda,h^{-1}]=\int_0^\lambda W^*(h'(y))\,dy,
\end{equation} 
where $y\in\Oma$ and $W^*(H)=HW(1/H)$. In the inverse-deformation formulation, we seek minimizers $h$ of \eqref{Es} that satisfy the boundary conditions
\begin{equation}\label{hbc}
	h(0)=0,\qquad h(\lambda)=1.
\end{equation}
We include the inequality constraint $h'\geq 0$. This ensures that the inverse strain is nonnegative and the orientation of matter is preserved. 

The inverse-deformation formulation is equivalent to the forward formulation when the deformation $f$ is bijective. However, the inverse-deformation formulation is better suited to the problem of fracture where $f$ can be discontinuous \cite{RHA, GH}. Cracks in a one-dimensional bar are characterized by jump discontinuities in the deformation map, viz. $F:=f'\sim \infty$, rendering both analysis and numerical computation difficult. In contrast, cracks appear as zero-slope regions in the inverse deformation, i.e. $H:=h'=0$. Even when the actual deformation $f$ may be discontinuous due to failure, the inverse deformation $h$ maintains continuity. Furthermore, the Lennard-Jones type potential $W$ appears as a double-well potential in $W^*$ with the two wells at $H=0$ and $H=1$, as shown in  Figure\,\ref{Wstar}. This is similar to phase-transition problems, but here in terms of the inverse strain \cite{RHA}. The $H=1$ well is the unbroken or `material phase' while the $H=0$ well represents empty space or the `vacuum phase'. Conservation of mass implies $H\rho=\rho^*$ where $\rho$ and $\rho^*$ are mass density in the reference and deformed configurations respectively. Consequently, the region in the deformed configuration where $H=0$ has zero mass. 

The energy in \eqref{Es} does not penalize the creation of new surfaces, i.e. jump discontinuities in $f$. The minimizers $h$ of \eqref{Es} may have arbitrary number of cracks and cannot sustain tension. This issue is resolved in \cite{RHA} by introducing a second-gradient energy in \eqref{Es}
\begin{equation}\label{energyse}
	E^*_\epsilon[\lambda,h]= \int_0^\lambda \frac{\epsilon^2}{2} (h'')^2+W^*(h')\,dy,
\end{equation}
where $\epsilon$ is a small parameter. It is shown in \cite{HPP} that \eqref{energyse} leads to Griffith's surface energy in the Gamma-limit as $\epsilon\searrow 0$. In any case, the problem reduces to solving for the minimizers of \eqref{energyse}, where $h$ satisfies the boundary conditions \eqref{hbc} and the unilateral constraint $h'\geq 0$. 

\subsection{Smoothness of weak solutions}\label{smooth}
We define the Hilbert space $\mathcal{H}:=H^2(\Oma,\Om)$. The admissible set is a closed convex cone defined as
\begin{equation}\label{defk}
	\mathcal{K}:=\{h\in\mathcal{H}: h(0)=0,\, h(1)=\lambda,\, h'\geq 0\}.
\end{equation}
Using the Sobolev embedding $\mathcal{H}\hookrightarrow C^1(0,\lambda)$, we characterize the inequality in \eqref{defk} pointwise. For any $h\in\mathcal{K}$, we define the broken set
\begin{equation}
	\mathcal{B}_h:=\{y\in[0,\lambda]: h'(y)=0\},
\end{equation} 
and the unbroken set 
\begin{equation}
	\mathcal{G}_h:=[0,\lambda]\backslash \mathcal{B}_h\equiv\{y\in[0,\lambda]: h'(y)>0\}.
\end{equation} 

Let $h,g\in\mathcal{K}$. The convexity of $\mathcal{K}$ implies that $(1-t)h+tg\in K$ for all $t\in[0,1]$. For a fixed $\lambda$, we define $e_\lambda(t):=E^*_\epsilon[\lambda,(1-t)h+tg]$. The function $h$ is an equilibrium solution if $[de_\lambda(t)/dt]|_{t=0}\geq0$. This gives the Euler-Lagrange inequality
\begin{equation}\label{ELineq}
	\int_0^\lambda \left\{\epsilon^2 h''(g''-h'')+S^*(h')(g'-h')\right \} \,dy\geq 0 \text{ for all } g\in \mathcal{K},
\end{equation}
where $S^*(H)=\frac{dW^*}{dH}.$ 

We derive the Euler-Lagrange equation in Propositions \ref{propEL} and \ref{propc4}. These results are established in \cite{RHA} in terms of the inverse strain. We formulate them here in terms of $h$ for completeness and to prove that the third derivative for a broken equilibrium is discontinuous, cf.\ Theorem \ref{jprop}.

\begin{proposition}\label{propEL}
	A point $h\in\mathcal{K}$ is a solution of \eqref{ELineq} if and only if there exist a nonegative Radon measure $\mu$ such that $\supp \mu\subseteq \mathcal{B}_h$ and $h,\mu$ satisfy
	\begin{equation}\label{EL}
		\int_0^\lambda \left\{\epsilon^2 h''\psi''+S^*(h')\psi'\right \} \,dy=\int_{\mathcal{B}_h}\psi'\,d\mu\text{ for all } \psi\in H^2(0,\lambda)\cap H^1_o(0,\lambda).
	\end{equation}
	Furthermore, $h\in C^2(0,\lambda)$.
\end{proposition}
\begin{proof}
Consider $\phi\in H^2(0,\lambda)\cap H^1_o(0,\lambda)$ such that $\phi'=0$ in $\mathcal{B}_h$. There exists sufficiently small $\alpha$ such that $g_1:=h+\alpha \phi$ and $g_2:=h-\alpha \phi$ are in $\mathcal{K}$. We substitute these into \eqref{ELineq} to find
\begin{equation} \label{mid}
		\int_0^\lambda \left\{ \epsilon^2 h''\phi''+S^*(h')\phi' \right\} \,dy= 0
	\end{equation}
for all such $\phi$. Let $\psi\in H^2(0,\lambda)\cap H^1_o(0,\lambda)$ such that $\psi'\geq 0$ in $\mathcal{B}_h$. We set $g=h+\psi$ in \eqref{ELineq} to obtain
	\begin{equation}\label{intbp}
		\int_0^\lambda \left\{\epsilon^2 h''\psi''+S^*(h')\psi'\right \} \,dy\geq 0
	\end{equation}
	for all such $\psi$. The left side of \eqref{intbp} is a linear functional in $\psi'\in C^0(0,\lambda)$.  By a version of the Riesz-Schwarz theorem \cite{kinder} and \eqref{mid}, there exist a nonegative Radon measure $\mu$ with $\supp \mu\subseteq \mathcal{B}_h$, such that $h,\mu$ satisfy \eqref{EL}.
	
	From \eqref{EL}, $h'''$ in the sense of distributions is given by
	\begin{equation}\label{th3dev}
		\epsilon^2 h'''=S^*(h')-C_1-\mu,
	\end{equation}
	where $C_1$ is a constant. Integrating \eqref{th3dev} gives
	\begin{equation}\label{hpp_pw}
		\epsilon^2 h''=\int_0^y S^*(h'(\xi))\, d\xi -C_1y-C_2-\nu,
	\end{equation}
	where $C_2$ is another constant and $\nu'=\mu$ in the distributional sense. Since $h\in C^1(0,\lambda)$ and $\nu$ is non-decreasing, $h''$ is continuous using the same argument as in \cite[Theorem 2.3]{RHA}.
\end{proof}

\begin{proposition}\label{propc4}
Let $h,\mu$ be a solution of \eqref{EL} and $G\subseteq \mathcal{G}_h$ be an open connected set. The solution $h$ restricted to $G$ is in $C^4(G,\mathbb{R})$ and satisfies
\begin{equation}\label{prop1}
	\epsilon^2 h^{iv}-\mathcal{M}^*(h')h''=0,
\end{equation} 	
where $\mathcal{M}^*:=\frac{d^2 W^* }{d H^2}$.
\end{proposition}
\begin{proof}
Since $\mu\equiv 0$ in $\mathcal{G}_h$, \eqref{th3dev} and \eqref{hpp_pw} imply
\begin{equation}\label{th3dev2}
		\epsilon^2 h'''=S^*(h')-C_1
	\end{equation}
in $G$. In view of $h\in C^2(0,\lambda)$, the right side of \eqref{th3dev2} is continuously differentiable. Differentiating \eqref{th3dev2} once gives \eqref{prop1}. 
\end{proof}

Although $h$ restricted to $\mathcal{G}_h$ is in $C^4$, we now show that $h'''$ is discontinuous at crack faces. A crack face $y_c\in[0,\lambda]$ separates the empty spaces from the material. It is defined such that $h'(y_c)=0$ and there exists no $\alpha\in \mathbb{R}$ such that $h'([y_c-\alpha,y_c+\alpha])=0$. Let $\mathcal{C}_h$ be the set of all crack faces of $h$. We assume no crack face is degenerate. For all $y_c\in\mathcal{C}_h$ there exists $\alpha>0$ such that $[y_c,y_c+\alpha)\in\mathcal{B}_h$ and $(y_c-\alpha,y_c)\in\mathcal{G}_h$ or vice-versa. This excludes the case when $h'=0$ at an isolated singular point and two crack faces touch each other. 

Let $\phi$ be a field in $\Oma$. We define the notation $\llbracket \phi \rrbracket$ for $y_c\in\mathcal{C}_h$ as
$$\llbracket \phi\rrbracket:=\lim_{y\rightarrow y_c^+}\phi-\lim_{y\rightarrow y_c^-}\phi.$$

\begin{theorem}\label{jprop}
	Let $h,\mu$ be a solution of \eqref{EL} and $y_c\in\mathcal{C}_h$ be a crack face. Then $h'''$ and $\mu$ are piecewise continuous with each having a jump discontinuity at $y_c$. Furthermore, at $y_c$ 
	\begin{equation}\label{jc}
		\epsilon^2\llbracket h'''\rrbracket=-\llbracket \mu \rrbracket\neq 0.
	\end{equation} 
\end{theorem}
\begin{proof}
	Consider $B\subseteq \mathcal{B}_h$ be open and connected. Since $h'=0$ in $\mathcal{B}_h$, we deduce $h'''=0$ in $B$. Without loss of generality, we assume that $y_c\in \mathcal{C}_h$ is a crack face such that for $\alpha>0$, $[y_c,y_c+\alpha)\in\mathcal{B}_h$ and $(y_c-\alpha,y_c)\in\mathcal{G}_h$. Using \eqref{th3dev2}, we deduce
	\begin{equation}
	\lim_{y\rightarrow y_c^-} \epsilon^2 h'''= S^*_0-C_1,\quad	\lim_{y\rightarrow y_c^+} \epsilon^2 h'''=0,
	\end{equation} 
	where $S^*_0=S^*(0)$ is a constant. Equation \eqref{hpp_pw} is satisfied pointwise. For $y\in B$, $h''(y)=0$. Substituting this in \eqref{hpp_pw} gives
	\begin{equation}\label{nu_temp}
		\nu=\int_0^y S^*(h'(\xi))\, d\xi -C_1y-C_2
	\end{equation}
	is differentiable within $B$. We differentiate \eqref{nu_temp} with respect to $y$ and get
	\begin{equation}
		\mu=S^*_0-C_1 \text{ in } B.
	\end{equation}
	We already know $\mu \equiv 0$ in $\mathcal{G}_h$. This implies,
	\begin{equation}
	\lim_{y\rightarrow y_c^-} \mu=0,\quad	\lim_{y\rightarrow y_c^+} \mu=S^*_0-C_1,
	\end{equation} 
and \eqref{jc}.

Assume that $h'''$ is continuous at $y_c$, i.e. $S^*_0= C_1$ and $\mu\equiv 0$. The solution $h\in C^4(0,\lambda)$. The inverse strain $H:=h'$ satisfies the equation
\begin{equation}\label{PLTeq}
	\epsilon^2 H''=S^*(H)-S_0,
\end{equation}
due to \eqref{th3dev2}, and the conditions
\begin{equation}\label{PLTbc}
	H(y_c)=0,\qquad H'(y_c)=0.
\end{equation}
The function $S^*(H)-S_0$ is continuously differentiable in $H$ because $W^*\in C^2([0,\infty)).$ Due to the Picard-Lindel\"of Theorem for second-order equations (\cite{CL}), there exists a unique function $\tilde{H}$ that satisfies \eqref{PLTeq} and \eqref{PLTbc} in an interval $I:=(y_c-\delta, y_c+\delta),\,\delta>0$. This function is $\tilde{H}\equiv 0$ by inspection. Since \eqref{PLTeq} is autonomous, we use the same argument with a new point $y\in I$ also satisfying \eqref{PLTbc} and confirm the uniqueness of this solution in larger domains. This gives $H\equiv 0$ is the unique solution of \eqref{PLTeq} and \eqref{PLTbc} in $[0,\lambda]$. Clearly, $\int_0^\lambda H\,dy=0$ and $h\notin \mathcal{K}$. This is a contradiction implying that $h'''$ is discontinuous.  
\end{proof}

\subsection{Work done by external forces}\label{workext}
The Euler-Lagrange equation \eqref{EL} embodies the configurational force balance \cite{eshelby1970, eshelby1975, gurtin2000}. In this section, we derive the actual force balance equation using the forward formulation. Our goal is to find an expression for work done on a part of the body by the forces external to it.

Assume for this subsection that $f:\Om\rightarrow\Oma$ is sufficiently smooth. By definition,
$$h=f^{-1},\qquad h'(y)=1/f'(h(y))$$
for all $y\in \Oma$. On differentiating further, we determine
\begin{equation}
	h''(y)=-\frac{f''(h(y))}{(f'(h(y)))^3}.
\end{equation}   
Substituting these in \eqref{energyse}, we find the total energy equals\footnote{See \cite{RHA, carlson}.} 
\begin{equation}
	E_\epsilon[\lambda,f]:=E^*_\epsilon[\lambda,f^{-1}]=\int_0^1 \frac{\epsilon^2}{2}\frac{f''^2}{f'^5}+W(f')\,dx.
\end{equation}
We determine the force balance equation by evaluating
\begin{equation}
	\left.\frac{d}{d\alpha}E[\lambda,f+\alpha \eta]\right|_{\alpha=0}=\int_0^1 \epsilon^2 \left\{ \frac{f''}{f'^5}\eta'' - \frac{5}{2}\left(\frac{f''}{f'^3}\right)^2\eta' \right\}+S(f')\eta'\,dx=0,
\end{equation}
where $\eta$ is an admissible variation in $f$ and $S:=\frac{dW}{dF}$. Using integration by parts, we find
\begin{equation}\label{lmb}
\begin{aligned}
	&\left[\epsilon^2 \frac{f''}{f'^5}\eta' \right]_0^1-\left[\left\{\epsilon^2 \left(\frac{f''}{f'^5}\right)'+ \frac{5}{2}\epsilon^2\left(\frac{f''}{f'^3}\right)^2-S\right\}\eta \right]_0^1+\\
	&\int_0^1\left\{\epsilon^2 \left(\frac{f''}{f'^5}\right)''+5\epsilon^2\left(\frac{f''}{f'^3}\right)\left(\frac{f''}{f'^3}\right)'-S' \right\}\eta \,dx =0.
\end{aligned}
\end{equation}
The integrand in the second line of \eqref{lmb} provides the strong form of the force balance equation. The first line of \eqref{lmb} equals the work done by external forces at the boundary. 

Consider a domain $[x_1,x_2]\subseteq\Om$. Let $\delta_1$ and $\delta_2$ be the change in displacements and $\upsilon_1$ and $\upsilon_2$ be the change in strain at points $x_1$ and $x_2$ respectively. The work done on the domain $[x_1,x_2]$ by external forces equals
\begin{equation}
\begin{aligned}\label{wd0}
	\mathcal{W}=&\left.\epsilon^2 \frac{f''}{f'^5}\right|_{x=x_2}\upsilon_2-\left.\left\{\epsilon^2 \left(\frac{f''}{f'^5}\right)'+ \frac{5}{2}\epsilon^2\left(\frac{f''}{f'^3}\right)^2-S\right\}\right|_{x=x_2}\delta_2\\
	&-\left.\epsilon^2 \frac{f''}{f'^5}\right|_{x=x_1}\upsilon_1+\left.\left\{\epsilon^2 \left(\frac{f''}{f'^5}\right)'+ \frac{5}{2}\epsilon^2\left(\frac{f''}{f'^3}\right)^2-S\right\}\right|_{x=x_1}\delta_1.
\end{aligned}
\end{equation}
In terms of the inverse deformation, \eqref{wd0} simplifies to
\begin{equation}\label{wd}
\begin{aligned}
	\mathcal{W}=&-\left.\epsilon^2 H'H^2\right|_{y=f(x_2)}\upsilon_2+\left.\left\{\epsilon^2\left(H''H-\frac{H'^2}{2}\right)+W^*-HS^*\right\}\right|_{y=f(x_2)}\delta_2\\
	&+\left.\epsilon^2 H'H^2\right|_{y=f(x_1)}\upsilon_1-\left.\left\{\epsilon^2\left(H''H-\frac{H'^2}{2}\right)+W^*-HS^*\right\}\right|_{y=f(x_1)}\delta_1,
\end{aligned}
\end{equation}
where $H=h'$. The Cauchy stress in the bar is
\begin{equation}\label{cauchy}
	\sigma=\epsilon^2\left(H''H-\frac{H'^2}{2}\right)+W^*-HS^*,
\end{equation}
and is independent of $y$ due to force balance \cite[Equation 5.6]{RHA}. Hence, \eqref{wd} simplifies to
\begin{equation}\label{wd2}
	\mathcal{W}=-\left.\epsilon^2 H'H^2\right|_{y=f(x_2)}\upsilon_2+\left.\epsilon^2 H'H^2\right|_{y=f(x_1)}\upsilon_1+\sigma(\delta_2-\delta_1).
\end{equation}
We use \eqref{wd2} in the next section to determine the entropy produced in the system during a quasistatic motion.

\section{Entropy production during quasistatic motion}\label{entropy}
\subsection{Dissipation and admissibility}
Consider a quasistatic, time-dependent family of equilibrium solutions $h(y,t): \Oma(t)\times [0,T]\rightarrow \Om$, where 
$$\Oma(t):=[0,\lambda(t)],\qquad h(0,t)=0,\qquad h(\lambda(t),t)=1.$$ 
The function $\lambda(t)$ gives the stretch value at time $t$. No inertia effects are considered here.

\begin{figure}[h]
  \centering
\includegraphics[width=0.45\textwidth]{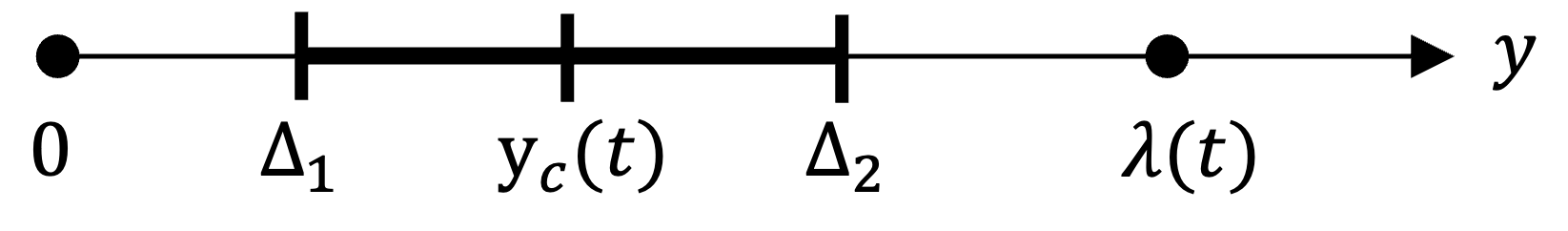}
  \caption{The domain $\mathcal{D}_*:=[\Delta_1,\Delta_2]\subset \Oma$ with a crack face at $y=y_c(t)$. }
  \label{quasi}
\end{figure} 
Let $\mathcal{D}_*:=[\Delta_1,\Delta_2]\subset \Oma(t)$ be a control volume as shown in Figure\,\ref{quasi}. Let $(')$ represent partial derivative with respect to $y$. The point $y_c(t)\in \mathcal{D}_*$ is a crack face such that $h'>0$ for $y\in[\Delta_1,y_c(t))$ and $h'=0$ for $y\in[y_c(t),\Delta_2]$ for all $t\in[0,T]$. The elastic energy of the material in $\mathcal{D}_*$ equals
\begin{equation}
	U^*_\epsilon= \int_{\Delta_1}^{\Delta_2} \frac{\epsilon^2}{2} (H')^2+W^*(H)\, dy.
\end{equation}
 where $H=h'$. 
 
 Let $(\dot{\,})$ represent partial derivative with respect to $t$. We use the Reynold's Transport Theorem and find the rate of change of internal energy \begin{equation}
 \begin{aligned}
 	\dot{U}^*_\epsilon=& \int_{\Delta_1}^{y_c(t)}\frac{\partial}{\partial t} \left \{\frac{\epsilon^2}{2} (H')^2+W^*(H)\right\}\, dy +\int^{\Delta_2}_{y_c(t)}\frac{\partial}{\partial t} \left \{\frac{\epsilon^2}{2} (H')^2+W^*(H)\right\}\, dy\\
 	&\quad +\left.\left(\frac{\epsilon^2}{2} (H')^2+W^*(H)\right)u\right|_{y=\Delta_1}^{y=\Delta_2}\\
 	=& \int_{\Delta_1}^{y_c(t)}\epsilon^2 H'v''+S^*(H)v'\, dy + \int^{\Delta_2}_{y_c(t)}\epsilon^2 H'v''+S^*(H)v'\, dy\\
 	&\quad +\left.\left(\frac{\epsilon^2}{2} (H')^2+W^*(H)\right)u\right|_{y=\Delta_1}^{y=\Delta_2},
 \end{aligned}
 \end{equation}
 where $v=\dot{h}$ and $u(y,t): \Oma(t)\times [0,T]\rightarrow \Om$ is the velocity of the material point at $y\in\Oma(t)$. We separate the integral because $H''$ has a jump discontinuity at $y=y_c(t)$, cf. Theorem \ref{jprop}. We use integration by parts and find
  \begin{equation}
  	\begin{aligned}
  		\dot{U}^*_\epsilon=& \int_{\Delta_1}^{y_c(t)}\left\{-\epsilon^2 H''+S^*(H)\right\}v'\, dy + \int^{\Delta_2}_{y_c(t)}\left\{-\epsilon^2 H''+S^*(H)\right\}v'\, dy\\
  		&+\epsilon^2 H'v'\bigg|_{y=\Delta_1}^{y=\Delta_2} +\left.\left(\frac{\epsilon^2}{2} (H')^2+W^*(H)\right)u\right|_{y=\Delta_1}^{y=\Delta_2}.
  	\end{aligned}
  \end{equation}
 Using integration by parts once again yields  
 \begin{equation}\label{ud}
  	\begin{aligned}
  		\dot{U}^*_\epsilon=& \llbracket \epsilon^2 H'' v\rrbracket +\int_{\Delta_1}^{y_c(t)}\{\epsilon^2 H''-{S^*(H)}\}'v\,dy+\int_{y_c(t)}^{\Delta_2}\{\epsilon^2 H''-{S^*(H)}\}'v\,dy \\
  		 	&+(-\epsilon^2H''v+S^*(H)v+\epsilon^2 H'v')\bigg|_{y=\Delta_1}^{y=\Delta_2} +\left.\left(\frac{\epsilon^2}{2} (H')^2+W^*(H)\right)u\right|_{y=\Delta_1}^{y=\Delta_2}.
  	\end{aligned}
  \end{equation}
The Euler-Langrange equation implies that $\{\epsilon^2 H''-S^*(H)\}'=0$ in $(\Delta_1(t),y_c(t))$. And $\{\epsilon^2 H''-S^*(H)\}'=0$ in $(y_c(t),\Delta_2(t))$ because here $H=0$. Consequently, the two integrals in \eqref{ud} are both zero.

The function $h$ maps $y\in\Oma$ to a material point $x\in\Om$  such that $x=h(y,t)$. We differentiate this equation with respect to time and get
\begin{equation}\label{fbc}
	0=Hu+v,\quad \implies v=-Hu.
\end{equation} 
The partial derivative of \eqref{fbc} with respect to $y$ gives
\begin{equation}\label{dbc}
v'=-H'u-Hu'.
\end{equation}
We substitute \eqref{fbc} and \eqref{dbc} into \eqref{ud} to find
\begin{equation}\label{ud2}
\begin{aligned}
	\dot{U}^*_\epsilon=& \llbracket \epsilon^2 H'' v\rrbracket-\epsilon^2 H'Hu'\bigg|_{y=\Delta_1}^{y=\Delta_2}+ \left.\left(\epsilon^2 H''H-\frac{\epsilon^2}{2} (H')^2+W^*(H)-HS^*(H)\right)u\right|_{y=\Delta_1}^{y=\Delta_2}.
\end{aligned}
\end{equation}
Equations \eqref{cauchy} and \eqref{ud2} give
\begin{equation}\label{ud3}
	\begin{aligned}
	\dot{U}^*_\epsilon=& \llbracket \epsilon^2 H'' v\rrbracket+\left(-\epsilon^2 H'Hu'+\sigma u\right)\bigg|_{y=\Delta_1}^{y=\Delta_2}.	
	\end{aligned}
\end{equation}

We now determine the instantaneous rate of work done on $\mathcal{D}_*$ from \eqref{wd2}.
This rate, denoted by $P$, equals 
\begin{equation}\label{p0}
	P=\left(-\epsilon^2 H'H^2\Upsilon+\sigma u\right)\bigg|_{y=\Delta_1}^{y=\Delta_2},
\end{equation}
where 
\begin{equation}\label{ras}
	\Upsilon:=\frac{du}{dx}
\end{equation}
is the rate of change of actual strain. We multiply \eqref{ras} with $H$ on both sides and get $H\Upsilon={u'}$. Substituting this into \eqref{p0} gives
\begin{equation}\label{p1}
	P= \left(-\epsilon^2 H'Hu'+\sigma u\right)\bigg|_{y=\Delta_1}^{y=\Delta_2}.
\end{equation}
The rate of energy dissipation is defined as
\begin{equation}\label{diss}
	D=P-\dot{U}^*_{\epsilon}=-\llbracket \epsilon^2 H'' v\rrbracket.
\end{equation}
Let $x_c(t)=h(y_c(t),t)$ be the location of the crack in the reference domain. The velocity of the crack in this domain is defined as
\begin{equation}\label{Vdef}
	V(t):=\dot{x}_c(t)=v(y_c(t),t)+H(y_c(t),t)\dot{y}_c (t).
\end{equation}
We substitute the expression of $v$ from \eqref{Vdef} in \eqref{diss} and get
\begin{equation}
	D=-\epsilon^2\llbracket H''\rrbracket V+ \epsilon^2\llbracket H''H\rrbracket\dot{y}_c.
\end{equation} 
The velocities $V$ and $\dot{y}_c$ do not depend on $y$ and were brought out of the double square brackets. From \eqref{cauchy}, we determine that $\llbracket H''H\rrbracket=0$. Hence,
\begin{equation}
	D=-\epsilon^2\llbracket H''\rrbracket V= \phi V.
\end{equation}  
where $\phi=-\epsilon^2\llbracket H''\rrbracket$ is the ``driving force" \cite{eshelby1951, eshelby1970, abey1990}. At a constant temperature, the dissipation rate $D$ is that temperature times the entropy production \cite{knowles79, podio2025assessment}. According to the second law of thermodynamics,
\begin{equation}\label{dissIneq}
	D=-\epsilon^2\llbracket H''\rrbracket V\geq 0.
\end{equation}
Similar  entropic inequality conditions were derived in problems of phase transitions where the deformation gradient is discontinuous \cite{abeyknowles88,abeyknowles92,abeyknowles93, abeyknowlesbook}. If given the sign of $\phi$ for an equilibrium $h$, the direction of $V$ is specified such that the inequality \eqref{dissIneq} is satisfied. Using \eqref{jc} and \eqref{dissIneq}, we determine that
\begin{equation}\label{ear}
	\phi=\llbracket \mu\rrbracket,\qquad D=\llbracket \mu\rrbracket V\geq 0.
\end{equation}
Since $\mu$ is positive in $\mathcal{B}_h$ and zero in $\mathcal{G}_h$, the sign of $\phi$ and subsequently of $V$ are easily obtained for a given $h$. In the following two subsections we discuss the cases of an end crack and an internal crack and determine the criterion of admissibility. These two cases cover all possible situations in the failure of an elastic bar. 

\subsection{End cracks}\label{endc}
\begin{figure}[h]
  \centering
\includegraphics[width=0.45\textwidth]{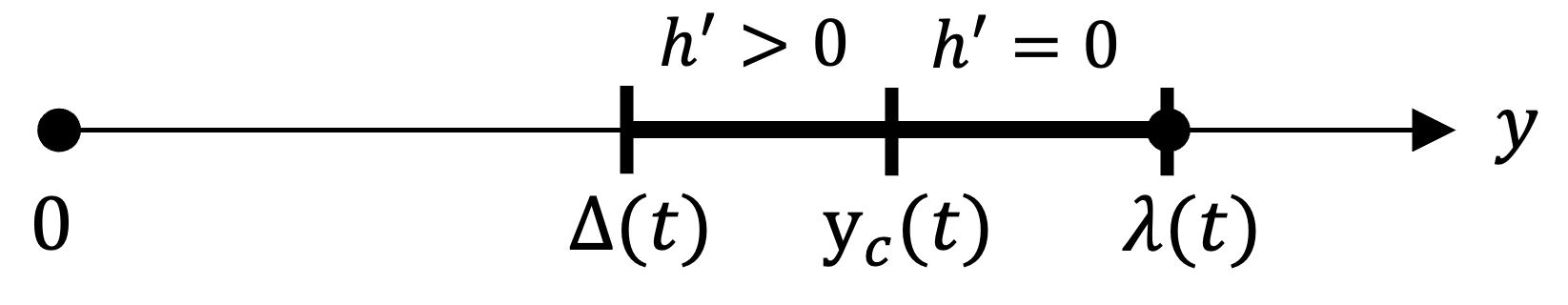}
  \caption{The domain $[\Delta_1(t),\lambda(t)]\subset \Oma$ with a crack face at $y=y_c(t)$. The picture shows the deformed configuration of an end crack located at $x_c=1$.}
  \label{endcrack}
\end{figure}

We consider the domain $[\Delta(t),\lambda(t)]\subset \Oma$, as shown in Figure\,\ref{endcrack}. The point $y_c(t)$ is the crack face of an end crack such that the region $[y_c(t),\lambda(t)]$ is empty space. That is, $h'(y)=0$ for all $y\in [y_c(t),\lambda(t)]$. Since $\mu$ is positive in $\mathcal{B}_h$ and zero otherwise, we establish that $\llbracket \mu\rrbracket>0$ at $y=y_c(t)$. From \eqref{ear}, $V\geq 0$ for entropic admissibility. 

Since $h(\lambda(t),t)=1$ and $h'(y)=0$ for all $y\in [y_c(t),\lambda(t)]$, the location of crack $x_c(t):=h(y_c(t),t)$ equals $1$. If $V>0$, then $x_c$ increases beyond $1$ and leaves the domain $\Om$. This is kinematically inadmissible. Hence,  $V=0$ necessarily. 

Similarly, if the end crack was at $x_c=0$, the sign of the driving force $\phi$ is switched and $V\leq 0$. However, $V<0$ is kinematically inadmissible using the same argument. We conclude once again that $V=0$.

\subsection{Internal cracks}\label{intc}
\begin{figure}[h]
  \centering
\includegraphics[width=0.45\textwidth]{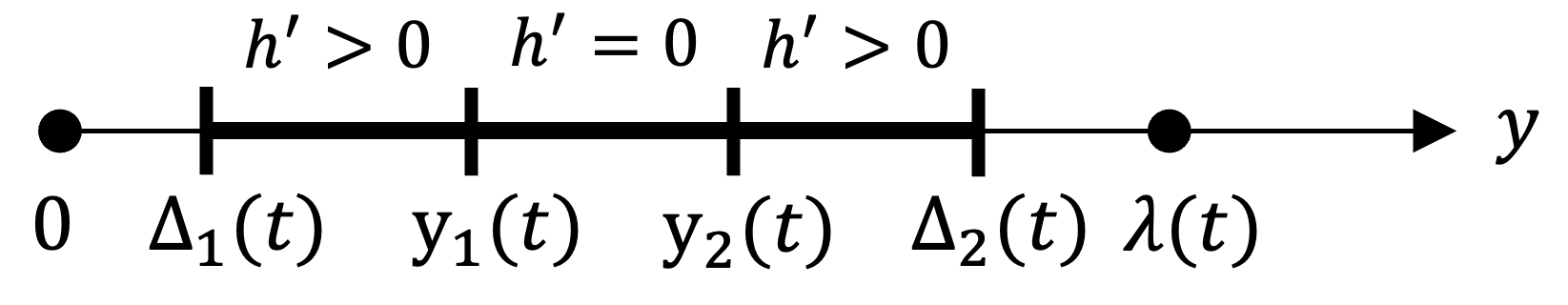}
  \caption{The domain $[\Delta_1(t),\lambda(t)]\subset \Oma$ with a crack face at $y=y_c(t)$. The picture shows the deformed configuration of an internal crack located at $x=x_c$.}
  \label{internal_crack}
\end{figure}

Consider the domain $[\Delta_1(t),\Delta_2(t)]\subset \Oma$, as shown in Figure\,\ref{internal_crack}. The points $y_1(t)$ and $y_2(t)$ represent the two crack faces of an internal crack such that $h'(y)=0$ for all $y\in[y_1(t),y_2(t)]$. Let $\phi_1=\left.\llbracket \mu\rrbracket\right|_{y=y_1(t)}$ and $\phi_2=\left.\llbracket \mu\rrbracket\right|_{y=y_2(t)}$ be the driving forces at the two crack faces. Using the sign of $\mu$, we determine that $\phi_1>0$ and $\phi_2<0$. From \eqref{ear}, the material velocities at the two crack faces satisfy $V_1\geq 0$ and $V_2\leq 0$.

Since $h'(y,t)=0$ for all $y\in[y_1(t),y_2(t)]$, the crack faces map to the same material point, i.e. $h(y_1(t),t)=x_c(t)=h(y_2(t),t)$. This implies, the two material velocities satisfy $V_1=\dot{x}_c(t)=V_2$. This is true if and only if $V_1=V_2=0$.

\subsection{Irreversibility condition}
The material velocity $V=\dot{x}_c$ describes the motion of the crack in the reference domain. The analyses in sections \ref{endc} and \ref{intc} show that $V=0$ for all cracks. This implies that if a crack nucleates at $x=x_c$ in the reference domain, this breakage point remains the same throughout the motion. New cracks may nucleate during a motion but the material locations of previously nucleated cracks do not change. We define the crack set as
\begin{equation}
	\Gamma(t):=\{x\in[0,1]: x=h(y,t) \text{ for } y\in\mathcal{B}_h\}.
\end{equation} 
The entropic inequality implies the \textit{irreversibility condition}
\begin{equation}\label{irre}
	\Gamma(s)\subseteq \Gamma(t) \text{ for all } s\leq t. 
\end{equation}
This condition is applicable to motions where the crack faces are separated from each other at all times. We assumed in section \ref{smooth} that the crack faces are not degenerate and $h'=0$ in a region of non-zero length. 

\section{Conditions for local stability}\label{condstb}
\subsection{Preliminaries for stability analysis}
A solution $h\in\mathcal{K}$ of \eqref{EL} is locally stable if it is a local minimizer, i.e. there exists an open neighborhood $U$ of $h$ such that $E^*_\epsilon[\lambda,g]\geq E^*_\epsilon[\lambda,h]$ for all $g\in U$. In this section, we prove second-order conditions for the stability of a quasistatic equilibria. The total energy $E^*_{\epsilon}: \mathcal{H}\rightarrow \mathbb{R}$ in \eqref{energyse} at a given value of $\lambda$ is twice Fr\'echet differentiable such that, 
\begin{equation}\label{sv}
	\delta^2 E_\epsilon^*[\lambda,h;\eta,\eta]=\int_0^\lambda \{\epsilon^2 (\eta'')^2+ \mathcal{M}^*(h')(\eta')^2\}\,dy \text{ for } \eta\in H^2(0,\lambda)\cap H^1_o(0,\lambda).
\end{equation} 
The definition of Fr\'echet differentiability, cf. \cite{vainberg}, gives 
\begin{equation}\label{frech1}
	E^*_{\epsilon}[\lambda,h+\alpha\eta]-E^*_{\epsilon}[\lambda,h]=\alpha\delta E^*_{\epsilon}[\lambda,h;\eta]+\alpha^2\delta^2E^*_{\epsilon}[\lambda,h;\eta,\eta]+o(||\alpha\eta||^2)
\end{equation}

Stability results are governed by what space of admissible variations is considered. We first define the admissible space in the absence of the irreversibility condition
\begin{equation}\label{vh1}
	\mathcal{V}_h^1=\{\eta\in H^2(0,\lambda)\cap H^1_o(0,\lambda): \eta'(y)\geq 0 \text{ for all } y\in\mathcal{B}_h\}	.
\end{equation}
The stability criterion using this admissible space is called Criterion 1.

We now define the space of admissible variations that satisfy the irreversibility condition. During a quasistatic motion, $H(y_c(t),t)=0$ for all crack faces $y_c\in\mathcal{C}_h$. We substitute this in \eqref{Vdef} and find
\begin{equation} \label{Vcond_irre}
	V=\left.\frac{\partial h}{\partial t}\right|_{y=y_c(t)}=0.
\end{equation}
The equality to zero is due to the irreversibility condition. This gives the space of admissible variations
\begin{equation}\label{admis}
	\mathcal{V}_h^2:=\{\eta\in H^2(0,\lambda)\cap H^1_o(0,\lambda): \eta'(y)\geq 0 \text{ for all } y\in\mathcal{B}_h,\, \eta(y_c)= 0 \text{ for all } y_c\in\mathcal{C}_h \}. 
\end{equation}
We study stability using this space in section \ref{c2} and call it Criterion 2.

\subsection{Criterion 1: Stability in the absence of irreversibility}\label{c1}
\begin{theorem}[Criterion 1]\label{c1th}
	Let $h\in\mathcal{K}$ be an equilibrium solution satisfying \eqref{EL}. Consider the space of admissible variations $\mathcal{V}_h^1$ defined in \eqref{vh1}. If $h$ is locally stable then $\delta^2 E_\epsilon^*[\lambda,h;\eta,\eta]\geq 0$ for all $\eta\in \mathcal{V}_h^1$ satisfying $\eta'=0$ in $\mathcal{B}_h$. Moreover, if $\delta^2 E_\epsilon^*[\lambda,h;\eta,\eta]> 0$ for all $\eta\in \mathcal{V}_h^1$ satisfying $\eta'=0$ in $\mathcal{B}_h$, then $h$ is a locally stable subject to Criterion 1.
\end{theorem}
\begin{proof}

$(\Rightarrow)$ Let $h$ be a local minimizer and $\eta\in \mathcal{V}_h^1$ such that $\eta'=0$ in $\mathcal{B}_h$. There exists $\alpha\in \mathbb{R}$ such that 
\begin{equation}\label{energyineq1}
	E^*_{\epsilon}[\lambda,h+\alpha\eta]\geq E^*_{\epsilon}[\lambda,h].
\end{equation}
Since $h$ satisfy \eqref{EL}, $\delta E^*_\epsilon[\lambda,h;\eta]=\int_{\mathcal{B}_h}\eta'\,d\mu=0$. Therefore, \eqref{energyineq1} and \eqref{frech1} imply $\delta^2 E_\epsilon^*[\lambda,h;\eta,\eta]\geq 0$.

$(\Leftarrow)$ Let $\delta^2 E_\epsilon^*[\lambda,h;\eta,\eta]> 0$ for all $\eta\in \mathcal{V}_h^1$ satisfying $\eta'=0$ in $\mathcal{B}_h$. For such $\eta$, $\delta E^*_\epsilon[\lambda,h;\eta]=\int_{\mathcal{B}_h}\eta'\,d\mu=0$. This implies, $E^*_{\epsilon}[\lambda,h+\alpha\eta]\geq E^*_{\epsilon}[\lambda,h]$ using \eqref{frech1} for sufficiently small $\alpha$. 

For variations $\eta \in \mathcal{V}_h^1$ with $\eta'>0$ in some $B\subseteq \mathcal{B}_h$, the first variation of the energy $\delta E_\epsilon^*[\lambda,h;\eta]=\int_{\mathcal{B}_h} \eta' d\mu>0$, cf.\ \eqref{EL}. Substituting this in \eqref{frech1} for sufficiently small $\alpha$ also gives $E^*_\epsilon[\lambda,h+\alpha\eta]> E^*_\epsilon[\lambda,h]$. Consequently, $E^*_\epsilon[\lambda,h+\alpha\eta]\geq E^*_\epsilon[\lambda,h]$ for all $\eta\in\mathcal{V}_h^1.$
\end{proof}

\begin{figure}[h]
  \centering
\includegraphics[width=0.4\textwidth]{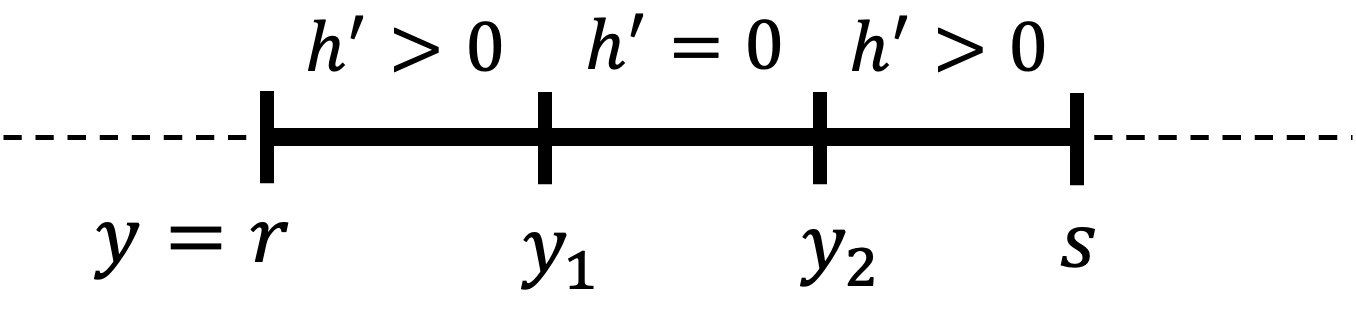}
  \caption{This figure depicts the location of points assumed in the proof of Proposition\ \ref{c1p1}. The points $y_1,\,y_2\in\mathcal{C}_h,\,i=1,2$ are crack faces and $r,\,s$ are points of maxima of $h'$.}
  \label{prop41}
\end{figure}

We use this result to prove the following proposition: 
\begin{proposition}\label{c1p1}
	Let $h\in\mathcal{K}$ be a broken equilibria satisfying \eqref{EL}. If $h$ has an internal crack, it is unstable subject to Criterion 1.
\end{proposition}  
\begin{proof}
We modify the method in \cite[Proposition 4.7]{RHA} to prove this result. Let $h\in \mathcal{K}$ be a broken solution on the $n$\textsuperscript{th} bifurcation branch, $n\geq 2$ that has an internal crack. The crack faces of this crack are $y_1,\, y_2\in\mathcal{C}_h$ as shown in Figure\ \ref{prop41}. Let $\lambda_*=\vol \mathcal{G}_h$ be the total length of all unbroken region. If the broken region is excised, the function $h'$ is even and $(2\lambda_*/n) $-periodic \cite{RHA}. At the boundary, $h''(0)=h''(\lambda)=0$. Let $r\in[0,y_1),$ and $s\in(y_2,\lambda]$ be such that $h''(r)=h''(s)=0$, $r,\,s \notin\mathcal{C}_h$. The points $r,\,s$ are points of maxima of $h'$ and exist on both sides of an internal crack due to the periodicity of the excised $h'$.  
	
	We define $c=h'(r)=h'(s)>0,$ and the function
	\begin{equation}
		\varphi(y)=\begin{cases}
			h'(y)-c, & y\in [r,s],\\
			0, & \text{otherwise.} 
		\end{cases}
	\end{equation}  
	Clearly $\varphi\in\mathcal{V}^1_h$ because $\varphi(0)=\varphi(\lambda)=0$ and $\varphi'=0$ in $\mathcal{B}_h$. Let $\psi \in \mathcal{V}^1_h$  such that $\psi'=0$ in $\mathcal{B}_h$, $\psi'(r)= 1$, and $\psi'(s)=0$. We define $\eta\in\mathcal{V}^1_h$ such that
	\begin{equation}
		\eta=\varphi-\tau\psi
	\end{equation}
	where $\tau$ is a small parameter. We evaluate \eqref{sv} for variation $\eta$ and get
	\begin{equation}
	\begin{aligned}
		\delta^2 E[\lambda, h; \eta,\eta]= &\int_r^{s} \{\epsilon^2 (h''')^2+ \mathcal{M}^*(h')(h'')^2\}\,dy \\
		&-2\tau \int_r^{s} \{\epsilon^2 h'''\psi''+ \mathcal{M}^*(h')h''\psi'\}\,dy + O(\tau^2).
	\end{aligned}
\end{equation}
The solution $h$ satisfy \eqref{prop1} in $\mathcal{G}_h$ and $h'=\psi'=0$ in $\mathcal{B}_h$. Using this with integration by parts gives
\begin{equation}\label{rhapf}
	\delta^2 E[\lambda, h; \eta,\eta]= 2\tau\epsilon^2 h'''(r)+O(\tau^2).
\end{equation}
We know that $h'''(r)<0$ because $r$ is a maxima of $h'$. Therefore, \eqref{rhapf} is negative for sufficiently small $\tau>0$. The solution $h$ is unstable due to Theorem\ \ref{c1th}.
\end{proof}

This proof does not work for solutions with only one or two end cracks\footnote{For a pictorial representation of these solutions, confer $1A$, $1B$ and $2B$ of Figure\, \ref{n1to4}.} because we are unable to find two separate maxima $r,\,s\in[0,\lambda]$ of $h'$. Unlike in \cite{carr}, $r,\,s$ cannot be the point of minima here when the solution is broken. In those cases, $r,\,s\in\mathcal{C}_h$. If $\psi'(r)\neq 0$, the first variation of the total energy is positive and the second variation of the energy is inconsequential for the stability in that direction, cf. Theorem\ \ref{c1th}.

These solutions with one or two end cracks remain the only candidates whose stability subject to Criterion 1 is yet to be determined. In Theorem\ \ref{c1th}, the second-order condition is checked for $\eta\in\mathcal{V}_h^1$ such that $\eta'=0$ in $\mathcal{B}_h$. Since these solutions only have end-cracks, the boundary conditions imply $\eta(y_c)=0$ for the crack faces $y_c\in \mathcal{C}_h$. This implies $\eta\in \mathcal{V}_h^2$. Furthermore, by definition $\mathcal{V}_h^2\subseteq \mathcal{V}_h^1$. The stability of these solutions using Criterion 1 and 2 is equivalent. We postpone their study to the next section where we present stability results for all equilibria using Criterion 2.

\subsection{Criterion 2: Stability in the presence of irreversibility}\label{c2}

\begin{theorem}[Criterion 2]\label{stabth}
	Let $h\in\mathcal{K}$ be an equilibrium solution satisfying \eqref{EL}. Consider the space of admissible variations $\mathcal{V}_h^2$ defined in \eqref{admis}. If $h$ is locally stable then $\delta^2 E_\epsilon^*[\lambda,h;\eta,\eta]\geq 0$ for all $\eta\in \mathcal{V}_h^2$ satisfying $\eta=0$ in $\mathcal{B}_h$. Moreover, if $\delta^2 E_\epsilon^*[\lambda,h;\eta,\eta]> 0$ for all $\eta\in \mathcal{V}_h^2$ satisfying $\eta=0$ in $\mathcal{B}_h$, then $h$ is a locally stable subject to Criterion 2.
\end{theorem}
\begin{proof}
$(\Rightarrow)$ Let $h$ be a local minimizer and $\eta\in \mathcal{V}_h^2$ such that $\eta=0$ in $\mathcal{B}_h$. There exists $\alpha\in \mathbb{R}$ such that 
\begin{equation}\label{energyineq2}
	E^*_{\epsilon}[\lambda,h+\alpha\eta]\geq E^*_{\epsilon}[\lambda,h].
\end{equation}
Since $h$ satisfy \eqref{EL}, $\delta E^*_\epsilon[\lambda,h;\eta]=\int_{\mathcal{B}_h}\eta'\,d\mu=0$. Therefore, \eqref{energyineq2} and \eqref{frech1} imply $\delta^2 E_\epsilon^*[\lambda,h;\eta,\eta]\geq 0$.

$(\Leftarrow)$	We first consider variations $\eta \in \mathcal{V}_h^2$ with $\eta'>0$ in some $B\subseteq \mathcal{B}_h$. Using \eqref{EL}, the first variation of the energy $\delta E_\epsilon^*[\lambda,h;\eta]=\int_{\mathcal{B}_h} \eta' d\mu>0$. Substituting this in \eqref{frech1} for sufficiently small $\alpha$ gives $E^*_\epsilon[\lambda,h+\alpha\eta]> E^*_\epsilon[\lambda,h]$ for all such $\eta$ irrespective of the sign of $\delta^2E^*_{\epsilon}[\lambda,h;\eta,\eta]$.
		
	The admissible variations $\eta\in\mathcal{V}_h^2$ that remain satisfy the condition $\eta'=0$ in $\mathcal{B}_h$. For these $\eta$, $\delta E_\epsilon^*[\lambda,h;\eta]=0$ and the sign of $\delta^2E^*_{\epsilon}[\lambda,h;\eta,\eta]$ determines stability. Since $\eta(y_c)= 0 \text{ for all } y_c\in\mathcal{C}_h $, we get $\eta=0$ in $\mathcal{B}_h$. Let $\delta^2 E_\epsilon^*[\lambda,h;\eta,\eta]>0$ for all such $\eta$. Using \eqref{frech1}, we get $E^*_\epsilon[\lambda,h+\alpha\eta]\geq E^*_\epsilon[\lambda,h]$ for sufficiently small $\alpha$. We conclude that the solution $h$ is locally stable. 
\end{proof}

Although the admissible variations are in $\mathcal{V}_h^2$, the second-order condition is satisfied by a minimizer in a smaller space of variations. We define this space as
\begin{equation}
	\mathcal{S}_h^2:=\{\eta\in H^2(0,\lambda)\cap H^1_o(0,\lambda):\, \eta=0 \text{ in } \mathcal{B}_h\}.
\end{equation}

Theorem\ \ref{stabth} has more physical implications. Cracks in a broken equilibrium $h$ may result in multiple unbroken parts that are separated by empty space from each other. Theorem \ref{stabth} implies that an equilibrium $h$ is stable if each of these unbroken parts are \textit{independently} stable. We show this in the following corollary.

\begin{corollary} \label{cor1}
	Let $h\in\mathcal{K}$ satisfy \eqref{EL}. Let ${G_i}:=(r_i,s_i), \,i=1,\dots, N$ be open, disjoint, connected subsets of $\mathcal{G}_h$ such that $r_i,\,s_i\in\mathcal{C}_h\cup\{0,\,\lambda\}$, and $\bigcup_{i=1}^N G_i=\mathcal{G}_h$. If $\delta^2 E_\epsilon^*[\lambda,h;\eta_i,\eta_i]>0$ for all $\eta_i\in \mathcal{S}_h^2$, $\supp \eta_i\subseteq G_i$  for all $i\in \{1,\dots,N\}$, then $h$ is locally stable.
\end{corollary}
\begin{proof}
Let ${G_i}:=(r_i,s_i)\subseteq \mathcal{G}_h$ as defined in the statement and $\eta\in\mathcal{S}_h^2$. The second variation in \eqref{sv} equals
\begin{align}\label{coreq1}
	\delta^2 E_\epsilon^*[\lambda,h;\eta,\eta]&=\int_{\mathcal{G}_h} \{\epsilon^2 (\eta'')^2+ \mathcal{M}^*(h')(\eta')^2\}\,dy +\int_{\mathcal{B}_h} \epsilon^2 (\eta'')^2\,dy,\\
	&\geq \sum_{i=1}^N\int_{G_i} \{\epsilon^2 (\eta'')^2+ \mathcal{M}^*(h')(\eta')^2\}\,dy.
\end{align}
Let $\eta_i\in \mathcal{S}^2_h$ such that $\supp{\eta_i}\subseteq G_i,\, \eta_i=\eta$ in $G_i$. The equation \eqref{coreq1} simplifies to
\begin{equation}
\begin{aligned}
	\delta^2 E_\epsilon^*[\lambda,h;\eta,\eta]&\geq \sum_{i=1}^N\int_{G_i} \{\epsilon^2 (\eta'')^2+ \mathcal{M}^*(h')(\eta')^2\}\,dy\\
	&= \sum_{i=1}^N\int_{G_i} \{\epsilon^2 (\eta_i'')^2+ \mathcal{M}^*(h')(\eta_i')^2\}\,dy\\
	&=\sum_{i=1}^N \delta^2 E_\epsilon^*[\lambda,h;\eta_i,\eta_i].
	\end{aligned}
\end{equation}
Let $\delta^2 E_\epsilon^*[\lambda,h;\eta_i,\eta_i]> 0$ for all $\eta_i\in \mathcal{S}^2_h$, $\supp \eta_i\subseteq G_i$ for all $i\in \{1,\dots,N\}$, then $\delta^2 E_\epsilon^*[\lambda,h;\eta,\eta]>0$ for all $\eta\in\mathcal{S}^2_h$. Theorem \ref{stabth} implies $h$ is locally stable. 
\end{proof}

The subsets $G_i$ defined in Corollary \ref{cor1} are the unbroken parts of $h$ separated by empty regions. If the second variation in \eqref{sv} is positive for all variations restricted to an unbroken part, we say that this  part is \textit{independently} stable. Corollary \ref{cor1} allows us to determine the stability of the full solution by studying the unbroken parts individually.

Consider the unbroken part that has crack at its two ends. Let $G:=(r,s)\subseteq\mathcal{G}_h$ such that $r,\,s\in\mathcal{C}_h$. Consider the function
\begin{equation}
	\xi:=\begin{cases}
		h', & y\in G,\\
		0, & y\in \Oma\backslash G.
	\end{cases}
\end{equation}
Since $h\in C^2(0,\lambda)$ and $r,\,s\in\mathcal{C}_h$, we conclude $\xi(r)=\xi'(r)=\xi(s)=\xi'(s)=0$. Therefore, $\xi\in\mathcal{S}^2_h$ with $\supp \eta\subseteq G$. For $\tau>0$, consider $\phi:=\xi+\tau\psi\in \mathcal{S}^2_h$ for some $\psi\in\mathcal{S}^2_h$. Evaluating \eqref{sv} for variation $\phi$ gives
\begin{equation}
\begin{aligned}
	\delta^2 E[\lambda, h; \phi,\phi]=&\int_0^\lambda \{\epsilon^2 (\xi'')^2+ \mathcal{M}^*(h')(\xi')^2\}\,dy + 2\tau \int_0^\lambda \{\epsilon^2 \xi''\psi''+ \mathcal{M}^*(h')\xi'\psi'\}\,dy \\
	& +\tau^2\int_0^\lambda \{\epsilon^2 (\psi'')^2+ \mathcal{M}^*(h')(\psi')^2\}\,dy,\\
	=& \int_r^s \{\epsilon^2 (h''')^2+ \mathcal{M}^*(h')(h'')^2\}\,dy + 2\tau \int_r^s \{\epsilon^2 h'''\psi''+ \mathcal{M}^*(h')h''\psi'\}\,dy \\
	&+ \tau^2\delta^2 E[\lambda, h; \psi,\psi].
\end{aligned}
\end{equation}
Restricted to $G$, the solution $h\in C^4(G,\mathbb{R})$, cf. Proposition \ref{propc4}. We use integration by parts in first two terms and get
\begin{equation}\label{svcal}
\begin{aligned}
	\delta^2 E[\lambda, h; \phi,\phi]=&  [\epsilon^2 h'''h'']_r^s+\int_r^s \{-\epsilon^2 h^{iv}h''+ \mathcal{M}^*(h')(h'')^2\}\,dy \\
	&+ 2\tau [\epsilon^2 h'''\psi']_r^s+ 2\tau \int_r^s \{-\epsilon^2 h^{iv}\psi'+ \mathcal{M}^*(h')h''\psi'\}\,dy \\
	&+ \tau^2\delta^2 E[\lambda, h; \psi,\psi].
\end{aligned}
\end{equation}
Since $r,\,s\in\mathcal{C}_h$ are crack faces, we know that $h''(r)=h''(s)=0$ and $\psi'(r)=\psi'(s)=0$. Moreover, $h$ satisfies \eqref{prop1} in $G$. Substituting these in \eqref{svcal} gives
\begin{equation}\label{svExi}
	\delta^2 E[\lambda, h; \xi+\tau\psi,\xi+\tau\psi]=\tau^2\delta^2 E[\lambda, h; \psi,\psi],
\end{equation} 
also implying that $\delta^2 E[\lambda, h; \xi,\xi]=0$. We show that $\xi$ is tangent to a one-parameter family of equilibria that have the same total energy \eqref{energyse}.
\begin{figure}[h]
  \centering
\includegraphics[width=0.4\textwidth]{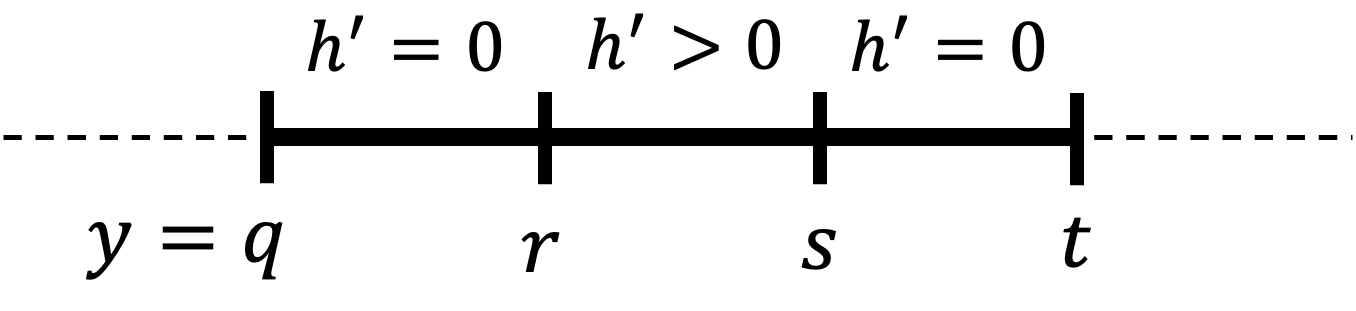}
  \caption{As defined in Theorem \ref{thfam}, the domain $G:=(r,s)\subseteq \mathcal{G}_h$ is unbroken region and region $R:=(q,t)$ such that $(q,r],\,[s,t)\in\mathcal{B}_h$. }
  \label{family}
\end{figure}
\begin{theorem}\label{thfam}
	Let $h$ be a solution of \eqref{EL} at stretch length $\lambda$ such that there exist an open, connected, unbroken region $G:=(r,s)\subseteq\mathcal{G}_h,$ $r,\,s\in\mathcal{C}_h$. Consider a region $R:=(q,t),\, q<r,\, t>s$ such that $R \backslash G\subseteq \mathcal{B}_h,\, h((q,r])=x_1,\, h([s,t))=x_2$ as shown in Figure\,\ref{family}. The one-parameter family of functions $T:(q-r, t-s)\rightarrow\mathcal{K}$ defined as
	\begin{equation}\label{fam}
	T(\theta)=
		\begin{cases}
			x_1, & y\in (q,r+\theta),\\
			h(y-\theta), & y\in (r+\theta, s+\theta),\\
			x_2,	& y\in (s+\theta, t),\\
			h(y), & y\in \Oma\backslash R,
		\end{cases}
			\end{equation}
	where $\theta\in (q-r, t-s)$, is a family of equilibrium solutions that have the same total energy.
\end{theorem}
\begin{proof}
Let $h_\theta:=T(\theta)$ for some $\theta\in (q-r, t-s)$. We first show that $h_\theta$ is an equilibrium solution. We evaluate \eqref{ELineq} for $h_\theta$ and get
	\begin{equation} \label{famEL}
	\begin{aligned}
		 \int_0^\lambda  &\left\{\epsilon^2 {h}''_\theta(g''-{h}''_\theta) +S^*({h}'_\theta)(g'-{h}'_\theta)\right \} \,dy =\\ & \int_{\Oma\backslash R} \left\{\epsilon^2 {h}''(g''-{h}'')+S^*({h}')(g'-{h}')\right \} \,dy \\
		  & + \int_{r+\theta}^{s+\theta} \left\{\epsilon^2 {h}''(y-\theta)(g''(y)-{h}''(y-\theta))+S^*({h}'(y-\theta))(g'(y)-{h}'(y-\theta))\right \} \,dy\\
		 & +S^*_0\int_q^{r+\theta}g'\,dy +S^*_0\int_{s+\theta}^t g'\,dy,
			\end{aligned}	
	\end{equation}
where $g\in \mathcal{K}$. We use the change of variable $z=y-\theta$ and define function $\varphi\in \mathcal{H}$, $\varphi(z)= g(z+\theta)$ in the second integral. Equation \eqref{famEL} becomes
\begin{equation}\label{famEL2}
	\begin{aligned}
		\int_0^\lambda  &\left\{\epsilon^2 {h}''_\theta(g''-{h}''_\theta) +S^*({h}'_\theta)(g'-{h}'_\theta)\right \} \,dy =\\ & \int_{\Oma\backslash R} \left\{\epsilon^2 {h}''(g''-{h}'')+S^*({h}')(g'-{h}')\right \} \,dy \\
		 & + \int_{r}^{s} \left\{\epsilon^2 {h}''(z)(\varphi''(z)-{h}''(z))+S^*({h}'(z))(\varphi'(z)-{h}'(z))\right \} \,dz\\
		 &+S^*_0(\varphi(r)-g(q))+S^*_0(g(t)-\varphi(s)),
	\end{aligned}
\end{equation} 
We construct $\zeta\in \mathcal{K}$ such that $\zeta=g$ in $\Oma\backslash R$ and $\zeta=\varphi$ in $G$. Substituting this in \eqref{famEL2} gives
\begin{equation}\label{famEL3}
	\begin{aligned}
		\int_0^\lambda  \left\{\epsilon^2 {h}''_\theta(g''-{h}''_\theta) +S^*({h}'_\theta)(g'-{h}'_\theta)\right \} \,dy =& \int_0^\lambda\left\{\epsilon^2 {h}''(\zeta''-{h}'')+S^*({h}')(\zeta'-{h}')\right \} \,dy\\
		\geq & 0.
	\end{aligned}
\end{equation} 
The inequality in \eqref{famEL3} is due to the fact that $h$ is an equilibrium solution, cf.\ \eqref{ELineq}. Therefore, $h_\theta$ are also equilibrium solutions. We find the total elastic energy \eqref{energyse} for $h_\theta$
\begin{equation}
\begin{aligned}
	E^*_\epsilon[\lambda,h_\theta]=& \int_{\Oma\backslash R} \left\{\frac{\epsilon^2}{2} (h'')^2+W^*(h')\right\}\,dy + \int_{r+\theta}^{s+\theta} \left\{\frac{\epsilon^2}{2} (h''(y-\theta))^2+W^*(h'(y-\theta))\right\}\,dy.
\end{aligned}
\end{equation}
This equals $E^*_\epsilon[\lambda,h]$ using a simple change of variable in the second integral.
\end{proof}
\begin{figure}[h]
  \centering
\includegraphics[width=\textwidth]{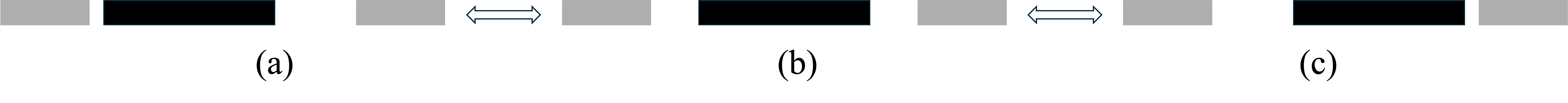}
  \caption{The black and grey rectangles represent the unbroken regions of the bar. The black region can translate within the empty space between the two grey regions.}
  \label{translate}
\end{figure}

It is not difficult to show that $T'(0)=-\xi$. The family of equilibria $T(\theta)$ defines the translational motion of the unbroken region $G$ within the empty space, as shown in Figure\,\ref{translate}. Such translations do not affect the solution energetically and are equivalent. An equilibrium $h$ is stable if \eqref{sv} is positive for all $\eta\in\mathcal{S}^2_h$ that are perpendicular to $\xi$ in the $L^2$ sense. We show this to prove the following sufficient condition:

\begin{theorem}[Sufficient condition using Criterion 2]\label{suffth}
	Let $h\in\mathcal{K}$ satisfy \eqref{EL}. Let $G_i:=(r_i,s_i)\subseteq\mathcal{G}_h,\, i=1,\,2,\,\dots,\,N$ be all unbroken regions such that $r_i,\,s_i\in\mathcal{C}_h$. Let $\xi_i\in\mathcal{S}^2_h$ be defined as
	\begin{equation}
	\xi_i:=\begin{cases}
		h', & y\in G_i,\\
		0, & y\in \Oma\backslash G_i.
	\end{cases}
\end{equation}
If $\delta^2 E_\epsilon^*[\lambda,h;\eta,\eta]>0$ for all $\eta\in \mathcal{S}^2_h$ such that $\int_{\Oma}\eta \xi_i\,dy=0$ for all $i=1,\,2,\,\dots,\,N$, then $h$ is locally stable.
\end{theorem}
\begin{proof}
	Consider the interval $G_i:=(r_i,s_i)$ defined in the theorem statement for some $i\in\{1,\,2,\,\dots,\,N\}$. Let $\langle\cdot,\cdot\rangle$ define the $L^2$ inner-product in $G_i$ such that $\langle u,v\rangle=\int_{r_i}^{s_i} uv\,dy$ for all $u,\,v\in L^2(G_i,\mathbb{R})$. We consider the variations $\varphi,\,\psi\in H^2_o(G_i, \mathbb{R})$ and note that
	\begin{equation}
		\delta^2 E_\epsilon^*[\lambda,h;\varphi,\psi]= \int_{G_i} \{\epsilon^2 \varphi''\psi''+ \mathcal{M}^*(h')\varphi'\psi'\}\,dy=\langle\varphi,A\psi\rangle=\langle A\varphi,\psi\rangle , 	
		\end{equation}
		where 
		\begin{equation}
			A\varphi=\epsilon^2 \varphi^{iv}-(\mathcal{M}^*(h')\varphi')' \text{ for } \varphi\in C^4_o(G_i,\mathbb{R})
		\end{equation}
		is a self-adjoint operator. 
		
		We consider the eigenvalue problem $A\phi=\beta\phi$. The eigenvalues $\beta_k,\, k=1,2,\dots$ are real and countable with no finite cluster point. The eigenfunctions $\phi_k$ are orthogonal in the $L^2$ sense and form a Schauder basis of $L^2(G_i,\mathbb{R})$. We point the interested reader to \cite[Chapter 7]{CL} for a detailed exposition of these results. 
		
		The function $\xi_i$ is an eigenfunction of $A$ with eigenvalue $0$, cf. \eqref{svExi}. Without loss of generality, assume that $\beta_1=0$ and $\phi_1=\xi_i/||\xi_i||_{L^2}$. Due to orthogonality of eigenfunctions, $\langle \xi_i,\phi_k\rangle =0$ for all $k\geq 2$. For $\psi\in H^2_0(G_i)$,
		\begin{equation}
			\psi=\sum_{k=1}^{\infty} \langle\psi,\phi_k\rangle\phi_k= a_1\xi_i+\zeta,
		\end{equation}
		where $a_1$ is a constant and $\zeta=\sum_{k=2}^{\infty} \langle\psi,\phi_k\rangle\phi_k$. Using \eqref{svExi}, we conclude 
		\begin{equation}
			\delta^2 E_\epsilon^*[\lambda,h;\psi,\psi]=\sum_{k=1}^{\infty} \langle\psi,\phi_k\rangle^2\delta^2 E_\epsilon^*[\lambda,h;\phi_k,\phi_k]=\sum_{k=2}^{\infty} \langle\psi,\phi_k\rangle^2\beta_k
		\end{equation}
		If $\delta^2 E_\epsilon^*[\lambda,h;\psi,\psi]>0$ for all $\psi\in H^2_o(G_i)$ such that $\langle\xi_i,\psi \rangle=0$ then the unbroken region $G_i$ is independently stable. If all unbroken regions are independently stable, the solution $h$ is locally stable due to Corollary \ref{cor1}.
\end{proof}

\section{Numerical stability results using Criterion 2} \label{num}
\subsection{Numerical implementation}
We solve \eqref{EL} using numerical continuation methods in the stretch parameter $\lambda$ employed in \cite{GH}. Since the domain $\Oma$ is dependent on $\lambda$, we change variables to solve on a fixed domain. We define $\Oms:=[0,1]$ and introduce the variable $s\in\Oms,\, s=y/\lambda$. Furthermore, we define $u\in H^2(0,1)\cap H^1_o(0,1),\, u'\geq -1$ such that $u(s):=h(\lambda s)-s$. The Euler-Lagrange equations \eqref{EL} gives
\begin{equation}\label{uEL}
	\int_0^1 \left\{\epsilon^2 u''\phi''+\lambda^3S^*\left(\frac{1+u'}{\lambda}\right)\phi'\right \} \,ds=\lambda^4\int_{\mathcal{B}_u}\{\mu\phi'+\zeta(u'+1)\}\,ds
\end{equation}    
for all $\phi\in H^2(0,1)\cap H^1_o(0,1)$, where $\mu$ is assumed to be in $L^\infty(0,1)$ and $\zeta\in L^\infty(0,1)$ is an admissible variation in $\mu$. The set $\mathcal{B}_u:=\{s\in\Oms: \, \lambda s\in \mathcal{B}_h\}$. The function $\mu$ now is the Lagrange multiplier for the inequality constraint $u'+1\geq 0$ and satisfies the Kuhn-Tucker conditions\footnote{See \cite{GH, NW}.} 
$$\mu\geq 0,\quad \mu(u'+1)=0.$$ 

It is clear that $u\equiv0$ is the homogeneous solution. A complete global bifurcation analysis of this problem is given in \cite{RHA}. Consider a solution $u$ on the $n$\textsuperscript{th}-bifurcation branch. If the broken region is excised, $u$ is odd and ${2\lambda_n^*}/{n\lambda}$ periodic, where $\lambda_n^*=\vol \mathcal{G}_h $ \cite{RHA}. Redistribution of the empty space $\mathcal{B}_h$ between crack faces gives equivalent solutions, cf. Theorem \ref{thfam}. Without loss of generality, we consider the case where the empty space is distributed to give $u$ odd and $2/n$ periodic. This allows us to compute the solution restricted to the domain $\Omn:=[0,1/n]$ and construct the complete solution $u$ on the domain $\Oms$ using reflection symmetry. 

Let the unknowns be $\bar{u}\in H^2(0,1/n)\cap H^1_o(0,1/n)$, $\bar{\mu}\in L^\infty(0,1/n)$ such that $\bar{u}=u$ and $\bar{\mu}=\mu$ in $\Omn$. We discretize the domain $\Omn$ in $N$ elements and use piecewise cubic Hermite polynomials to interpolate $\bar{u}$. This interpolation allows the approximation of $\bar{u}$ to be in $H^2(0,1/n)$. At each node, there are three unknowns: $\bar{u}_k,\,\bar{u}'_k,\,\bar{\mu}_k$, where subscript $k$ denotes the value of the function at the $k$\textsuperscript{th}-node. The right side of \eqref{uEL} is evaluated as in \cite{GH}:
\begin{equation}
	\lambda^4\int_{\mathcal{B}_{\bar{u}}}\{\bar{\mu}\bar{\phi}'+\bar{\zeta}(\bar{u}'+1)\}\,ds\approx \lambda^4\sum_{ s_k\in\mathcal{B}_{\bar{u}}^d} \{\bar{\mu}_k\bar{\phi}_k'+\bar{\zeta}_k(\bar{u}'_k+1)\}
\end{equation}
where $\mathcal{B}_{\bar{u}}^d:=\{s_k:\,k=1,\,2,\dots,\,N+1;\, \lambda s_k\in \mathcal{B}_{u}\}$ is the finite dimensional approximation of $\mathcal{B}_u$. The broken set $\mathcal{B}_{\bar{u}}^d$ is also unknown. We determine it along with $\bar{u}$ using the active-set method \cite[Section 4.1]{GH}. Once we obtain $\bar{u}$, we construct its $2/n$-periodic, odd extension and find $u$ in the domain $\Oms$. The function $\mu$ is obtained from the even extension of $\bar{\mu}$. 

\subsection{Stability analysis}
Given a solution $u$ in the $n$\textsuperscript{th}-branch, we determine its stability numerically using Theorem \ref{suffth}. We define the set of variations for $u$ equivalent to those for $h$ in $\mathcal{S}^2_h$ as
\begin{equation}
	\mathcal{S}^2_u:=\{v\in H^2(0,1)\cap H^1_o(0,1):\, v(s)=0 \text{ for all } s\in\mathcal{B}_u\}.
\end{equation}
The second variation of the total energy \eqref{sv} in terms of $u$ equals
\begin{equation}\label{usv}
	SV[\lambda,u;v,v]=\frac{1}{\lambda^3}\int_0^1 \left\{\epsilon^2 (v'')^2+ \lambda^2\mathcal{M}^*\left(\frac{1+u'}{\lambda}\right)(v')^2\right\}\,ds \text{ for } v\in \mathcal{S}_u.
\end{equation}

We consider the discretization of $\Oms=[0,1]$ with $nN$ elements. The solution $u$ obtained through our numerical procedure is in the form of a vector $U$ of length $2(nN+1)$. This vector is
\begin{equation}
	U=\begin{bmatrix}
		u_1 & u'_1 & u_2 & u'_2 & \dots & u_k & u'_k & \dots & u_{nN+1} & u'_{nN+1}
	\end{bmatrix}^T,
\end{equation}
where subscript $k$ denotes the value of the function at the $k$\textsuperscript{th}-node. We use the same interpolation for $v\in \mathcal{S}_u$ as we do for $u$. Therefore, $v$ is approximated using the vector
\begin{equation}
	V=\begin{bmatrix}
		v_1 & v'_1 & v_2 & v'_2 & \dots & v_k & v'_k & \dots & v_{nN+1} & v'_{nN+1}
	\end{bmatrix}^T,
\end{equation}
The values $v_k=v'_k=0$ for all $s_k\in\mathcal{B}^d_u$ by definition and $v_1=v_{nN+1}=0$ due to the boundary condition. Let $\tilde{V}$ be a vector of size $M\leq 2(nN+1)$ obtained by removing the known entries from vector $V$. The second variation in \eqref{usv} is approximated as
\begin{equation}
	SV[\lambda,u;v,v]\approx \tilde{V}^T\,G(U)\,\tilde{V}
\end{equation}
where $G(U)$ is a $M\times M$ matrix assembled using the finite-element discretization. The eigenvalues of $G(U)$ determine the stability of the solution. If the solution $u$ has $P$ number of connected, unbroken regions such that both ends of the region are crack faces, then $G(U)$ has at least $P$ zero eigenvalues, cf.\ Theorem \ref{suffth}. If the remaining eigenvalues of $G$ are positive, we conclude using the same theorem that $u$ is stable subject to Criterion 2. 

\subsection{Results}\label{stabresults}

\begin{figure}[h]
	\centering
	\begin{subfigure}[c]{0.48\textwidth}
		\centering
         \includegraphics[width=\textwidth]{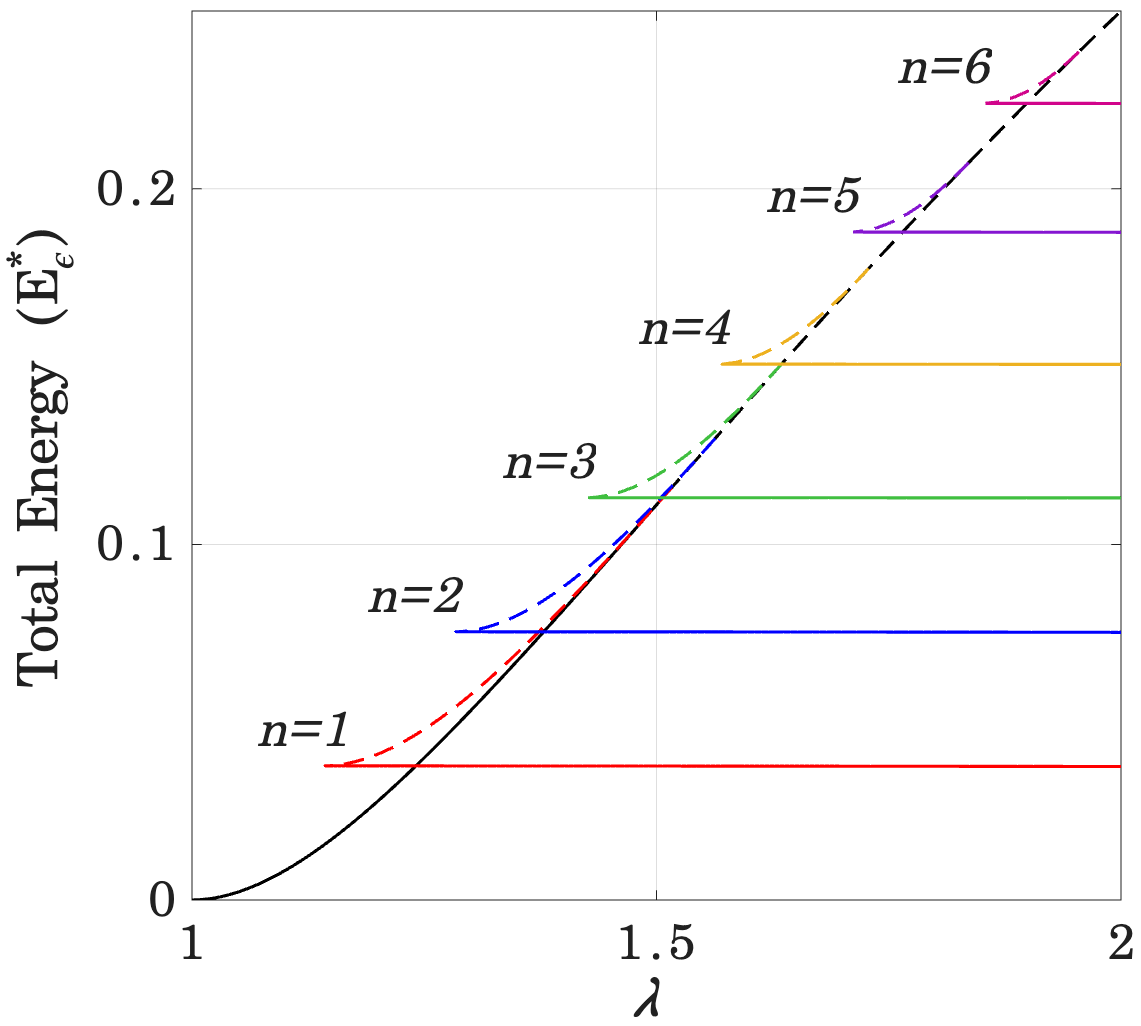}
         \caption{}
          \label{energy}
	\end{subfigure}
	\hfill
	\begin{subfigure}[c]{0.48\textwidth}
		\centering
         \includegraphics[width=\textwidth]{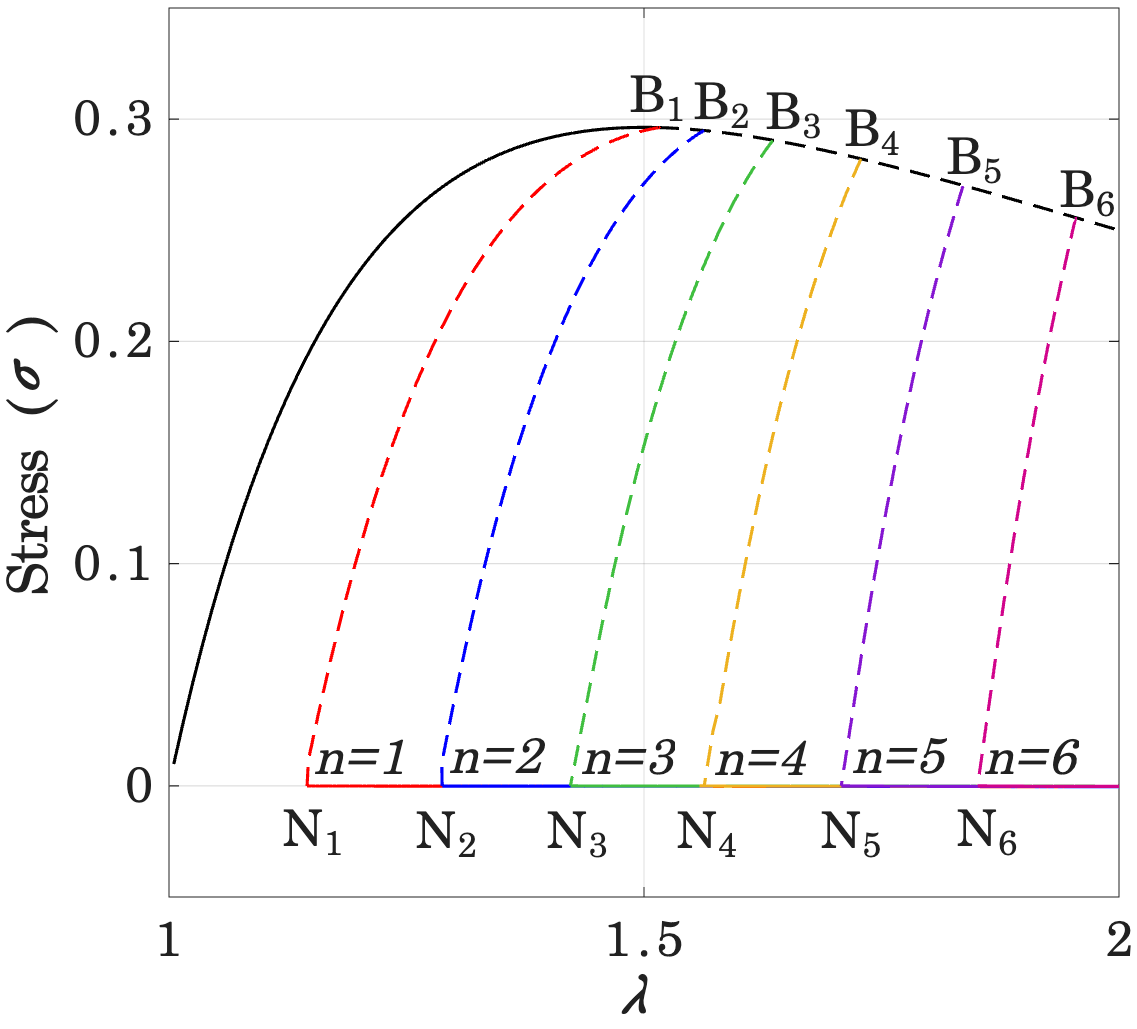}
         \caption{}
          \label{stress}
	\end{subfigure}
 \caption{(a) Total energy $E^*_\epsilon$ versus stretch length $\lambda$. (b) Stress $\sigma$ versus stretch length $\lambda$. The solid lines are locally stable and dashed lines are unstable solutions}
  \label{Fig1}
\end{figure}
 
We compute results for $\epsilon=0.1$ and the constitutive law
\begin{equation}
	W(F)=\left(1-\frac{1}{F}\right)^2,\quad \implies W^*(H)=H(1-H)^2.
\end{equation}
We compute the $n=1,\,2,\,\dots,\,6$ bifurcating branches. We use $N=600/n$ elements and relative and absolute tolerances of $10^{-9}$ in the Newton's method. Figure \ref{energy} shows the total energy $E^*_\epsilon$ and Figure\,\ref{stress} shows the stress $\sigma$ in the bar as a function of stretch length $\lambda$. The solid lines show stable solutions and the dashed lines show unstable solutions. The homogeneous solution $h\equiv y/\lambda,\,u\equiv 0$ is shown in black. It loses stability at critical load $\lambda=1.5163$. Subcritical pitchfork bifurcations branch out from the homogeneous solution and are shown in different colors. The two sides of the pitchforks are equivalent in terms of energy and stress and overlap each other in the two bifurcation diagrams. 

In Figure\,\ref{stress}, the bifurcation points on the black curve are marked with letters B\textsubscript{n}. Along the bifurcating branches, the value of stress decreases until it reaches zero at the point N\textsubscript{n}. These nonhomogeneous solutions with non-zero stress are not broken ($\mathcal{B}_h=\emptyset$) and are unstable. Cracks appear at a finite number of material points at N\textsubscript{n}. Beyond this, the solutions are broken ($\mathcal{B}_h\neq\emptyset$) and locally stable.  Since the stress is zero in all the broken solutions, the braches overlap each other in Figure\,\ref{stress} even though the solutions are different, cf. Figure\,\ref{energy}. 

We incorporate Theorem \ref{suffth} numerically and find that broken solutions in all bifurcation branches are locally stable subject to Criterion 2. Beyond $\lambda=1.2409$, the broken solutions in the branch $n=1$ have the lowest total energy and are energetically favorable. The broken solutions in $n\geq 2$ branches are local minimizers and have higher total energy. This is because more crack surfaces contain more surface energy. An unbroken bar undergoing incremental loading transitions from the homogeneous solution branch to a branch with broken solutions. This may happen at loadings at which the branch has lower energy than the homogeneous solution or when the homogeneous solution loses local stability. The total energy in the system drops and stress goes to zero. Due to imperfections, inertial effects, randomness or other real-world factors, the bar may transition into a broken solution of $n\geq 2$ branch and lock into a local minimizer. In those cases, internal or multiple cracks appear in the bar. 

\begin{figure}[H]
  \centering
\includegraphics[width=\textwidth]{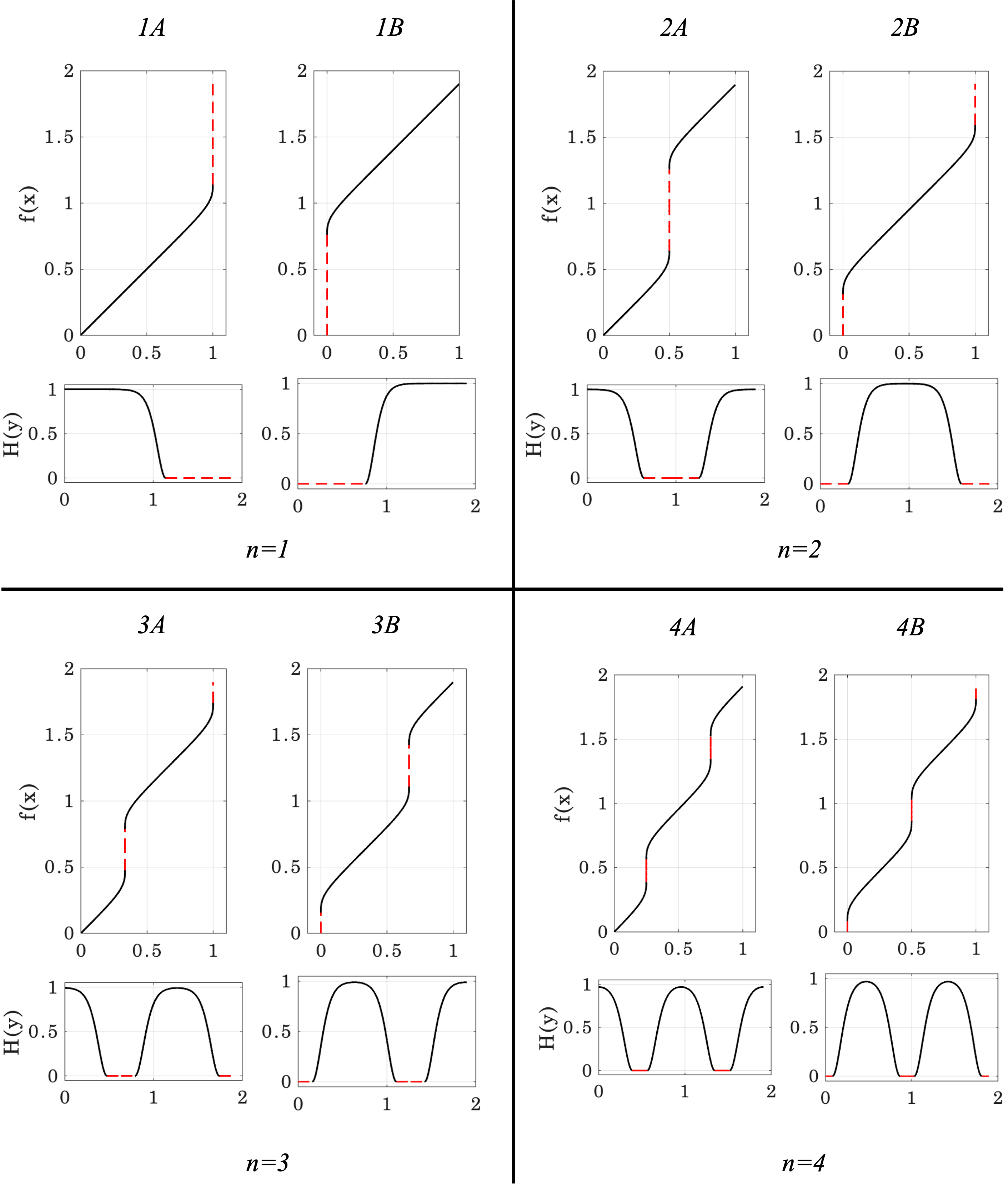}
  \caption{The deformation map $f$ and inverse strain $H$ for $n=1$ to 4 branches at $\lambda=1.9$. The letters $A$ and $B$ denote the two different sides of the same pitchfork. The black solid lines shows the material and the red dashed lines show empty space}
  \label{n1to4}
\end{figure}

\begin{figure}[H]
  \centering
\includegraphics[width=\textwidth]{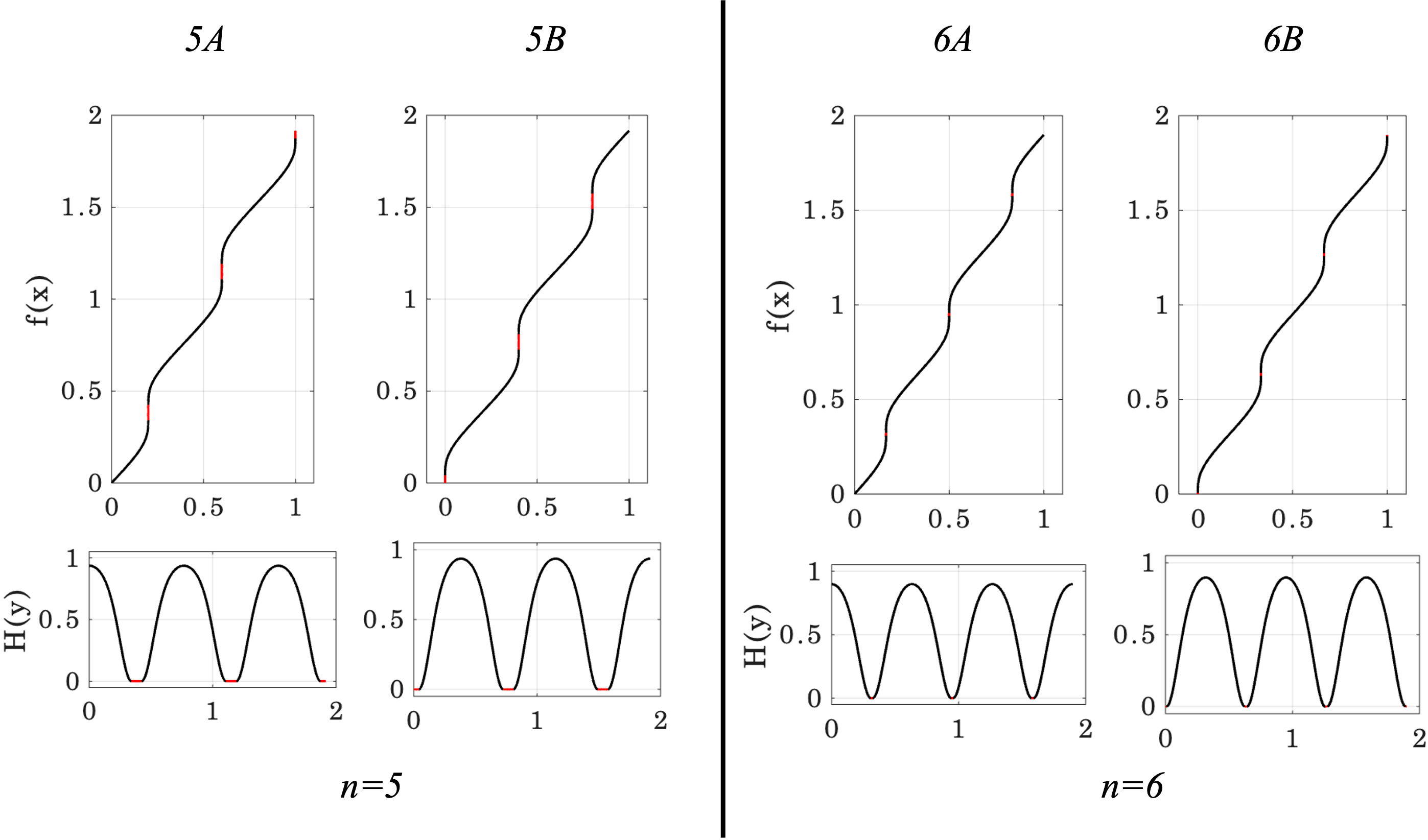}
 \caption{The deformation map $f$ and inverse strain $H$ for $n=5$ and 6 branches at $\lambda=1.9$. The letters $A$ and $B$ denote the two different sides of the same pitchfork. The black solid lines shows the material and the red dashed lines show empty space}
  \label{n5to6}
\end{figure}

We show the deformation map $f$ and inverse-strain $H$ of broken solutions at $\lambda=1.9$ for $n=1,\,2,\,\dots,\,6$ in Figures \ref{n1to4} and \ref{n5to6}. The letters $A$ and $B$ in the figures are used to identify the two different sides of the pitchfork branch. The black solid lines shows the material, i.e. $h'>0$, and the red dashed lines show empty space with $h'=0$. The crack faces precisely delineate the two regions.

In $1A$, the crack appears at the right end of the bar, i.e. at $x=1$. After the crack nucleates and $\lambda$ is increased, the hard loading device moves further away. The unbroken region remains as it is. In contrast, the crack is at $x=0$ in $1B$. The unbroken region is attached with the hard loading device and moves along with it as $\lambda$ is increased. The elastic, unbroken region in $1A,\,1B$ and $2A$ are attached with one of the boundaries. None of the broken pieces in these solutions can translate within the empty space. Therefore, in these cases $P=0$ and we find numerically that $G(U)$ has all positive eigenvalues. The solution $2B$ has two end cracks. The unbroken region in this case can translate within the empty space, cf. Theorem \ref{thfam}. In this case, $P=1$ and $G(U)$ has one zero and all other positive eigenvalues. The remaining solutions can be seen as a combination of the above. Solutions $3A,\,3B$ and $4A$ have one unbroken region each that can translate within the empty space. Solutions $4B,\,5A,\,5B$ and $6A$ have two such regions and $6B$ has three such unbroken regions. The matrix $G(U)$ in each of these cases have the same number of zero eigenvalues as the number of such unbroken regions.

Transitions between local minima on different branches must satisfy the irreversibility condition \eqref{irre}. Consequently, transitions to a lower $n$-valued branch from a higher one are not possible unless when the crack faces may come in contact with each other. The transitions that satisfy the irreversibility condition are towards a higher $n$-valued branch and require an increase in the total energy of the system. Such transitions are not possible without an external intervention that results in this increase. A broken bar undergoing incremental loading using a hard loading device remains locked into the same $n\geq 1$ branch. The hard-loading device continues to move forward and the empty space between the crack faces increases without any change in the total energy. 

\section{Concluding remarks} \label{conc}
We derive an irreversibility condition for the failure of one-dimensional, homogeneous, elastic bars as a consequence of the second law of thermodynamics. We study the stability of broken equilibria of the bar in the presence of this condition. The variations satisfy the entropic inequality, and the stability analysis models the real-world problem. We show that all broken equilibria of the elastic bar, including those with internal or multiple cracks, are locally stable.

The irreversibility condition found here was not known at the time of the publication of \cite{GH}. A similar analysis resulting in the same irreversibility condition is also possible for the composite bar studied in that work. The pseudo-rigid foundation introduces additional terms with lower order derivatives in the energy functional. These terms do not qualitatively impact the jump-discontinuity in the third derivative of the inverse deformation. The broken solutions presented in \cite{GH} are from the first critical mode and are stable in the absence of irreversibility. The presence of the irreversibility condition further limits the space of admissible variations and does not change the stability results reported in that work.  

Second law of thermodynamics in fracture assuming prescribed configurational forces is studied in \cite{gurtin1996, gurtin1998, gurtin2000}. An expression of the energy release rate available at the crack tip is derived in these works that is akin to the well known J-integral \cite{atkinson1968, rice1968J}. In the one-dimensional setting considered here, the bar breaks completely and the phenomena of crack propagation is not observed. This is not true for fracture of a two-dimensional solid employing the inverse-deformation approach. In that case, the second law will result in the expression for energy-release rate in addition to an irreversibility condition. We will pursue this question in a future work.

\section*{Acknowledgements}
I express my sincere gratitude to Timothy J. Healey and Phoebus Rosakis for stimulating and helpful conversations. I thank David Bindel, Anurag Gupta, and Sankalp Tiwari for helpful comments on a previous version of this article.

\bibliographystyle{siam}
\bibliography{thesisbib.bib}

\end{document}